\documentclass{amsart}

\newtheorem{theorem}{Theorem}[section]
\newtheorem{conjecture}{Conjecture}[section]
\newtheorem{lemma}[theorem]{Lemma}
\newtheorem{proposition}[theorem]{Proposition}
\newtheorem{corollary}[theorem]{Corollary}
\newtheorem{definition}[theorem]{Definition}
\newtheorem{example}[theorem]{Example}
\newtheorem{remark}[theorem]{Remark}

\newcommand \Set {\underline{Set}}
\newcommand \U {\mathcal{U}}
\newcommand \B {\mathcal{B}}
\newcommand \e {\epsilon}
\newcommand \g {\mathfrak{g}}
\newcommand \M {\mathcal{M}}
\newcommand \T {\mathbb{T}}
\newcommand \C {\mathbb{C}}
\def\pd#1#2{\dfrac{\partial#1}{\partial#2}}
\input xy
\xyoption{all}

\begin{document}

\title{ Hopf algebras with trace and representations}

\author{C. De Concini}\address{Dip. Mat. Castelnuovo, 
Univ. di Roma La Sapienza, Rome, Italy}\email{corrado@mat.uniroma1.it}
\author{C. Procesi}\address{Dip. Mat. Castelnuovo, 
Univ. di Roma La Sapienza, Rome, Italy}\email{claudio@mat.uniroma1.it} 
\author{N.Reshetikhin}\address{Department of
Mathematics, University of California, Berkeley, CA 94720 \\ USA}
\email{reshetik@math.berkeley.edu} 
\author{M.Rosso}\address{D.M.A., Ecole Normale Superieure 45,
rue d'Uim, 75230 Paris Cedex 05, France}\email{rosso@math.ens.fr}

\begin{abstract} We study the restriction of representations
of Cayley-Hamilton algebras to subalgebras. This theory is applied
to determine tensor products and branching rules for representations 
of quantum groups at roots of 1.

\end{abstract}
\maketitle 
\tableofcontents

\section*{Introduction}
%\addcontentsline{toc}{section}{Introduction}

Irreducible representations of quantized universal enveloping algebra
were classified in \cite{D1}{D2}. These algebras are finite  dimensional
over their center and are Cayley-Hamilton algebras \cite{P2}.

Here we study the restriction representations of Cayley-Hamilton
algebras to subalgebras. This theory is applied to the tensor product
of two generic representations of $U_\epsilon(\g)$, and of 
$U_\epsilon(\mathfrak{b})$, and to the branching of generic
irreducible representations when $U_\epsilon(\mathfrak{b})\subset
U_\epsilon(\g)$. here $\epsilon$ is an primitive root of 1 of an odd
degree.

The center of  $U_\epsilon(\g)$ has the central Hopf subalgebra 
$Z_0$ generated by $\ell$-th powers of root generators of and 
by $\ell$-th powers of generators of the Cartan subalgebra \cite{D1}.
This central Hopf subalgebra $Z_0$ is isomorphic to the algebra
of polynomial functions on $G^*$ which is a Poisson Lie dual to $G$.
Let $\pi: Spec(Z)\to Spec(Z_0)=G^*$ be the natural projection
induced by the inclusion of $Z_0$ to the center $Z$ of
$U_\epsilon(\g)$. According to the general theory of Cayley-Hamilton
algebras there exist a Zariski open subvariety $S\subset Spec(Z_0)$
such that $\pi$ is a finite covering map and the algebra is semisimple
over $S$. In case of $U_\epsilon(\g)$
this projection has $\ell^r$ fibers where $r$ is the rank of the 
Lie algebra $\mathfrak{g}$. Thus, each central character $\chi\in
Spec(Z)$ with $\pi(\chi)\in S$ defines an irreducible representation
$V_\chi$.

The tensor product $V_\chi\otimes V_{\chi'}$ is completely
reducible if the product $\pi(\chi)\pi(\chi')$ ( in $G^*$)
is generic, i.e. belongs to $S$. Our results imply:
\begin{equation}
V_\chi\otimes V_{\chi'}\simeq \oplus_{\chi''\in Spec(Z), / /
\pi(\chi'')=\pi(\chi)\pi(\chi')} (V_{\chi''})^{\oplus m}
\end{equation}
where $m=\ell^{|\Delta_+|-r}$. Here $|\Delta_+|$ is the
number of positive roots.

Similar decompositions hold for the restriction of $V_\chi$ to 
$U_\epsilon(\mathfrak{b})$ and for the tensor product of 
$U_\epsilon(\mathfrak{b})$-modules.

In sections 1 to 4 we recall the general theory of semisimple representations
of Cayley-Hamilton algebras (CH-algebras for short). Then in sections
5 to 7 we study how a semisimple
representation of a CH-algebra restricts to a CH-subalgebra. Then
we apply this theory to the decomposition of tensor product of
semisimple representations of the special class of Hopf algebras
that we call CH-Hopf algebras. Finally we study examples of such
Hopf algebras which are quantum groups at roots of 1. The special
form of the multiplicities suggests that for
many natural varieties related to the corresponding
Poisson Lie groups the Poisson tensor is constant in certain
birational coordinate system. Such coordinates are known in many cases.

The work of N.R. was supported by NSF grant DMS-0070931, C.P. and
N.R. would like to thank V. Toledano for interesting discussions.

\section{ $n-$dimensional representations. }

In this paper by {\bf ring} we mean  an associative ring with 1.
An algebra over a commutative ring $A$ will be  an associative and
unital algebra.
\smallskip

In this section we remind some basic facts of universal algebra.
We will use the following categories (we denote by $A$ a commutative ring):
\begin{itemize}
\item $\underline{Set}$ is the category of sets,
\item $C$ and $C(A)$ are categories of commutative rings and
commutative $A$-algebras
respectively,
\item $N$ and $N(A)$ are categories of non commutative rings and of
non commutative $A$-algebras respectively.
\end{itemize}

\subsection{The universal $n$-dimensional representation}
Given a ring $B$ one denotes by $M_n(B)$ the full ring of $n\times n$
matrices over $B$.
    If $f:B\to C$ is a ring homomorphism we can construct
$M_n(f):M_n(B)\to M_n(C)$ the homomorphism induced on matrices. We will use
systematically the
following simple well known Lemma.
\begin{lemma} If $I$ is a 2-sided ideal in $M_n(B)$ then
$I=M_n(J)$ for a (unique)
2-sided ideal
$J$ in $B$ (and $M_n(B)/I=M_n(B/J)$).
\end{lemma}

Furthermore, if $f:B\to C$ is such that $M_n(f):M_n(B)\to M_n(C)$
annihilates $I$, then there is a morphism
$\bar{f}:B/J\to C$ such that the following diagram commutes
\[
\begin{array}{lcl}
M_n(B) &\stackrel{f}\longrightarrow & M_n(C)\\
\quad\searrow & &\nearrow_{M_n(\bar{f})}   \\&
M_n(B/J)
\end{array}
\]
Let $R$ be a ring, by an $n-$dimensional representation of $R$ over
a commutative ring $B$ one means a homomorphism

\[
\phi: R\to M_n(B)
\]

We will use notation $\mathcal R^n_R(B):=\hom_N(R,M_n(B)) $ for the
set of all these
representations. It is clear that it defines a
functor $\mathcal R^n_R: {\mathcal C}\to \Set$. The image of a given ring
homomorphism $f:B\to C$ under this functor is the map ${\mathcal R}^n_R(B)\to
{\mathcal R}^n_R(C)$ obtained by composing with
\[
M_n(B)\stackrel{M_n(f)}\longrightarrow M_n(C).
\]

\begin{lemma} The functor $\mathcal R^n_R$ is representable.
\end{lemma}
\proof  We should prove that there is a commutative ring $A_n(R)$ such that
${\mathcal R}^n_R(B)=Hom_C(A_n(R),B)$.

Let $a_\alpha$ be a set of generators for the ring
$R$. This gives a presentation of it as a
quotient ring of the free (non commutative) algebra
${\mathbb Z}\langle x_\alpha \rangle$:
\[\pi:{\mathbb Z}
\langle x_\alpha \rangle\to R={\mathbb Z}\langle x_\alpha
\rangle/K,\ K:=Ker \pi.
\]

with generators $a_\alpha$ being images of $x_\alpha$.

For each
$\alpha$ choose a set of
$n^2$ variables $\xi_{i,j}^\alpha$. Let $A:={\mathbb
Z}[\xi_{i,j}^\alpha]$ be the
polynomial ring in
all these variables.

Define the {\it generic matrices}
$\xi_\alpha$ in $M_n(A)$ by setting $\xi_\alpha$ to be the matrix which,
in the $i,j$ entry, has
coefficient  $\xi_{i,j}^\alpha$.

Let $j:{\mathbb Z}\langle x_\alpha \rangle\to M_n(A)$ be the algebra
homomorphism defined by $j(  x_\alpha)=\xi_\alpha.$  Let finally $I$ be the
2-sided ideal in $M_n(A)$ generated by $j(K)$. By the previous Lemma,
$I=M_n(J)$
for some ideal $J$ of $A$ and thus we have the mapping $j_R:R\to M_n(A/J)$
and the commutative diagram
\[
\begin{array}{lcl}
{\mathbb Z}\langle x_\alpha \rangle&\stackrel{j}\longrightarrow & M_n(A)\\
\downarrow _\pi& & \downarrow _{M_n(p)}\\
R &\stackrel{j_R}\longrightarrow & M_n(A/J)
\end{array}
\]
Here $p:A\to A/J$ is the quotient map.

Let $\phi:R\to M_n(B)$ be any representation
$\phi(a_\alpha)=(f_{i,j}^\alpha)$. It gives the homomorphism
of commutative rings $A\to B$, $\xi_{i,j}^\alpha\mapsto
f_{i,j}^\alpha$ . This homomorphism induces the
required map $\overline{\phi}:A/J\to B$
for which  the diagram
\[
\begin{array}{lcl}
R &\stackrel{j_R}\longrightarrow & M_n(A/J)\\
   {_\phi} \searrow  & &\swarrow_{M_n(\overline{\phi})}\\
     &\ M_n(B)
\end{array}
\]
commutes. This proves representability of $\mathcal R^n_R$.

\qed

Notice that in the proof we constructed the commutative ring $A_n(R) := A/J.$

Notice that, for a finitely generated ring $R$ also $A_n(R)$ is finitely
generated. This construction is functorial.
   Indeed, for each ring
homomorphism
$f: R\to S$ we have a corresponding homomorphism of commutative rings
$A_n(f): A_n(R)\to A_n(S)$ defined naturally. It is clear that
the diagram
\[
\begin{array}{lcl}
R &\stackrel{j_R}\longrightarrow & M_n(A_n(R))\\
\downarrow_f &  &\downarrow _{M_n(A_n(f))}\\
S&\stackrel{j_S}\longrightarrow & {M_n(A_n(S))}
\end{array}
\]
is commutative.

\begin{definition}
We shall denote by $A_n(R)$ the universal ring $A/J$ and shall refer
to the mapping
\begin{equation}\label{j_R}
R\stackrel{j_R}\longrightarrow M_n(A_n(R))
\end{equation}
    as the
{\it universal $n-$dimensional representation}.
\end{definition}

Universality of the ring $A_n(R)$ means the commutativity of
the diagram above.

\begin{remark} The commutative ring $A_n(R)$ may be zero.
This means that the ring $R$ does not have any $n$-dimensional
representation.
\end{remark}
   Instead of working with rings we can work with algebras over a
commutative ring
     $A$.
Clearly all the discussion from above carries over.
Moreover if
$R$ is a finitely generated algebra so is the universal ring $A_n(R)$.
The universality of $A_n(R)$ implies the following theorem

\begin{theorem} The functor $B\to M_n(B)$ from the category
${\mathcal C}(A)$ of commutative $A-$algebras to
the category ${\mathcal N}(A)$ of non commutative  $A-$algebras has a
right adjoint
$$hom_{\mathcal N}(R,M_n(B))=hom_{\mathcal C}(A_n(R),B).$$
for each $R$ in $\mathcal N$.
\end{theorem}

\begin{example} Consider the ring $U$ generated by three
elements $H,X,Y$ with defining relations
\begin{equation}
HX-XH=Y
\end{equation}

\end{example}

In this case the ring $A_n(U)$ is the polynomial ring in the $2n^2$ variables
$h_{ij},x_{ij}$.
\begin{example} Consider the commutative polynomial ring  $\mathbb
Z[x,y]$  generated by two
elements $x$ and $y$.

\end{example}

In this case the ring $A_n(\mathbb Z[x,y])$ is generated by the
$2n^2$ variables
$x_{ij},y_{ij}$ modulo the quadratic equations $\sum_sx_{is}y_{sj}-
\sum_sy_{is}x_{sj}$. It is not known if these equations generate in
general a prime ideal!

\begin{example} Consider the $\mathbb Q$ algebra $U$ generated by two
elements $ X,Y$ with defining relations
\begin{equation}
   XY-YX =1
\end{equation}

\end{example}

In this case the ring $A_n(U)$ is 0. If instead we work over $\mathbb
Z$ we get not trivial rings, since the matrix
equation
$XY-YX =1 $  can be solved in some characteristic $p>0$.

\begin{example} Let $k$ be a field and $R:=M_m(k)$ the algebra of
$m\times m$ matrices.  Then $A_n(R)=0$ unless $m$ divides $n$;  in
the case $n=mr$, let $\overline k$ be an algebraic closure of $k$.
Consider the embedding $GL(r,\overline k)$ into
   $GL(mr,\overline k)$ given by the tensor product $ 1\otimes A:
\overline k^m\otimes \overline k^r\to \overline
k^m\otimes \overline k^r$. $A_n(R)$ is the coordinate algebra of the
homogeneous space   $GL(mr,\overline
k)/GL(r,\overline k)$ (cf. \cite{L1},\ \cite{L2}) \end{example}

    \vskip20pt

\subsection{Equivalence between universal $n$-dimensional representations}
When one studies representations
one has a natural equivalence given by {\it changing the basis}.
\footnote{In the
arithmetic theory one needs a
little care, one has to define the automorphism group of matrices as {\it
group scheme}.}   We shall not describe this theory for general
rings (cf. \cite{P}), but assume
now that all rings are algebras over a field
$k$.  Let $B$ be a $k$-algebra. An invertible matrix $g\in GL(n,k)$ defines
a $B-$automorphism by conjugation:
\[
M_n(B)\stackrel{C(g)}\longrightarrow M_n(B),\qquad C(g)(A):=gAg^{-1}
\]

Let $R$ be an algebra over $k$. From the universal property of
the universal $n$-dimensional representation, every matrix $g\in
GL(n,k)$ defines a homomorphism, $\overline g:A_n(R)\to A_n(R)$
making the following
    diagram  commutative
\[
\begin{array}{lcl}
R &\stackrel{j_R}\longrightarrow &M_n(A_n(R))\\
\downarrow j_R &  &\downarrow M_n (\overline g )\\
M_n(A_n(R))&\stackrel{C(g)}\longrightarrow& M_n(A_n(R))
\end{array}
\]
Notice that such $\overline g$ is unique due to the universality
of $A_n$.

Composing with some other $h\in G$ we get:
$$\begin{aligned}
    M_n( \overline  {hg}))\circ j_R &=C(h g)\circ j_R=C(h)\circ
C(g)\circ j_R= \\
C(h)\circ M_n( \overline g )\circ j_R &=M_n( \overline g )\circ
C(h)\circ j_R=
M_n( \overline g )\circ  M_n( \overline
h )\circ j_R =M_n({\overline g}\circ{\overline h})\circ j_R
\end{aligned}
$$
and therefore
\[
M_n(( \overline  {hg}))= M_n( {\overline g}\circ{\overline h} )
\]
which implies ${\overline {hg}}={\overline g}\circ{\overline h}$.
In other words we get an action of $GL(n,k)$ on $A_n(R)$
$g:a\mapsto \overline g^{-1}(a)$, an action of $GL(n,k)$ on $M_n(k)$
by conjugation
and the diagonal action $m\otimes_k a\mapsto C(g)m\otimes_k
\overline{g^{-1}}(a)$
on $M_n(A_n(R))=M_n(k)\otimes_k A_n(R).$

    The identity   $C(g)\circ j_R=M_n( \overline g)\circ j_R$  means that
$M_n( \overline g)^{-1}C(g)\circ j_R=  j_R$ or that
$j_R$ maps
$R$ into the elements which are $ GL(n,k)$ invariant:
\[
R\stackrel{j_R}\longrightarrow M_n(A_n(R))^{GL(n,k)}.
\]
One of the aims of the theory is to understand better the previous
map; it is clear
that in the algebra of invariants we
find all the {\it characters}  $Tr(j_R(a)),\ a\in R$ or even all coefficients
of characteristic polynomials. This
justifies introducing such characters formally as will be done in the next
section. In the meantime let us interpret the
construction for the free algebra. So let us assume that $k$ is an
infinite field which allows us to identify formal polynomials with functions.

For the free (non commutative) algebra $k\langle x_\alpha
\rangle_{\alpha\in I} $,
we have seen that for each
$\alpha$ we choose a set of
$n^2$ variables $x_{i,j}^\alpha$ and let $A_{n,I}:=k[x_{i,j}^\alpha]$ be the
polynomial ring in
all these variables.\smallskip

   We have seen that the universal map is the map
$j(  x_\alpha)=\xi_\alpha $ sending  each variable
$x_\alpha$ to the corresponding  generic matrix
$\xi_\alpha$,   the matrix which, in the $i,j$ entry, has
value  $x_{i,j}^\alpha$.\smallskip

The ring $A_{n,I}$ is best thought of as the ring $k[M_n(k)^I]$ of
polynomial functions on  $M_n(k)^I$ and the ring $M_n(A_{n,I} )$ is best
thought of as the ring of polynomial maps $f:M_n(k)^I\to M_n(k)$.

Now assume that $I$ is a set with $m$ elements. Choose it to be
$I=\{1,\dots, m\}$.
Then
\[
M_n(k)^I=M_n(k)^m=\{(\xi_1,\dots,\xi_m),\ \xi_i\in M_n(k)\}
\]
Define
\[
A_{n,m}=k[M_n(k)^m]=k[\xi_i^{hk}]
\]
The  generic matrix
$\xi_i$   is thus a {\it coordinate function} mapping $(\xi_1,\dots,\xi_m)\to
\xi_i$.   The $GL(n,k)$ action is
identified to the obvious action on functions:
\[
    (gf)(\xi_1,\dots,\xi_m):=gf(g^{-1} \xi_1 g,\dots,g^{-1}\xi_m g)g^{-1}
\]
Finally the ring  $M_n(A_{n,m} )^{GL(n,k)}$ can be identified to the
{\it ring of
${GL(n,k)}$ equivariant maps}:
\[
f:M_n(k)^n\to M_n(k)\,|\quad  f(g \xi_1 g^{-1},\dots,g \xi_m g^{-1})=
gf(\xi_1,\dots,\xi_m)g^{-1}.
\]

Let us denote by $C_n(\xi_1,\dots,\xi_m):=M_n(A_{n,m})^{GL(n,k)}$ the ring of
equivariant maps.

\section{Cayley-Hamilton algebras.}
\subsection{Algebras with trace. }   We start with a formal definition which
belongs to universal algebra.

\begin{definition}  An associative algebra with trace, over a commutative ring
$A$ is an associative algebra $R$ with a 1-ary operation
\[
t:R\to R
\]
which is assumed to satisfy the
following axioms:
\begin{enumerate}
\item  $t$ is $A-$linear.
\item  $t(a)b=b\,t(a), \quad\forall a,b\in R.$
\item  $t(ab)=t(ba), \quad\forall a,b\in R.$
\item $t(t(a)b)=t( a)t(b), \quad\forall a,b\in R.$
\end{enumerate}

This operation is called a formal trace.
We denote $t(R):=\{t(a),\ a\in R\}$ the image of $t$.
\end{definition}

\begin{remark} We have the following implications:

     Axiom 1) implies that $t(R)$ is an $A-$submodule.

     Axiom 2) implies that $t(R)$ is in the center of $R$.

     Axiom 3) implies that $t$ is 0 on the space of commutators $[R,R]$.

     Axiom 4) implies that $t(R)$ is an $A-$subalgebra and that $t$ is also
$t(R)-$linear.
\end{remark}

The basic example of algebra with trace is of course the algebra of $n\times n$
matrices over a commutative ring $B$ with the usual trace.

Notice that we have made no special requirements on the value of
$t(1)$.\bigskip

Algebras with trace form a category, where objects are algebras
with trace  and morphisms algebra homomorphisms which commute with trace
mappings. An ideal in a trace algebra is a {\it trace} ideal i.e. one
which is stable under the trace. Then the usual
homomorphism theorems are valid.

The subalgebra $t(R)$ is called the {\it trace algebra}.
\bigskip

\begin{example}
As usual in universal algebra the category of algebras with trace has
free algebras.
Given a set $I$ it is easily seen that the free
algebra with trace in the variables $x_i,\ i\in I$, is obtained as follows:
\end{example}
First one constructs the free algebra $A\langle x_i\rangle$ in the variables
$x_i,\ i\in I$ whose basis over $A$ are the free monomials in
these variables. Next one defines {\it cyclic equivalence} of
monomials where when we
decompose a monomial $M=AB$ we set $AB\cong BA$.

Finally for every equivalence class of monomials we pick a commutative variable
$t(M)$ and form finally the polynomial algebra:
\[
F\langle x_i\rangle:= A\langle x_i\rangle[t(M)],\ i\in I,
\]
where $M$ a monomial up to cyclic equivalence.

We set $T\langle x_i\rangle:= A[t(M)]$ the commutative polynomial algebra in
the variables $t(M)$. The trace is defined as the unique
$T\langle x_i\rangle-$linear map for which $t:M\to t(M)$.  It is
easily verified
that this gives the free algebra with trace.\smallskip
\begin{remark} There is an obvious base change construction on
algebras with trace. If $R$ is an $A-$ algebra with
trace and $B$ a commutative $A$ algebra then $R\otimes_AB$ has the
natural trace $t(r\otimes b):=t(r)b$. We have
$t(R\otimes_AB)=t(R)\otimes_AB$.
\end{remark}\bigskip

\subsection{$n$-dimensional representations of algebras with traces}

We can now apply the theory of \S 1 to algebras with traces. The only
difference
is now that, if $R$ is an algebra with  trace, by an
$n-$dimensional representation over
a commutative ring $B$ one means a homomorphism
$\phi: R\to M_n(B) $  which is compatible with traces, where matrices have the
standard trace.  The discussion of \S 1 can be repeated
verbatim, we have again a representable functor, a universal (trace preserving)
map $R\stackrel{i_R}\longrightarrow M_n(B_n(R))$, an action of $GL(n,k)$ on
$B_n(R)$ and a map
\[
R\stackrel{i_R}\longrightarrow M_n(B_n(R))^{GL(n,k)}
\]
We shall call $B_n(R)$ the {\it coordinate ring of the $n-$dimensional
representations of $R$} and the map $i_R$ the {\it generic
$n-$dimensional representation  of $R$.} \medskip

We shall see that, under suitable assumptions, this map $i_R$ is an
isomorphism. For this we first point
out some elements which are always in the kernel
of $i_R$.

\subsection{Cayley Hamilton algebras. }

      At this point we will restrict  the discussion to the case in
which $A$ is a
field of characteristic 0. The positive characteristic
theory can to some extent be developed, provided we start the axioms
from the idea
of a {\it norm} and not a trace. Since the theory is
still incomplete we will not go into it now.\bigskip

The basic algebraic restriction which we know for the algebra of
$n\times n$ matrices over a commutative ring $B$ is the Cayley Hamilton
theorem:

Every matrix $M$ satisfies its characteristic polynomial
$\chi_M(t):=\det(t-M)$.\smallskip

The main remark that allows to pass to the formal  theory is that,
in characteristic 0, there are universal polynomials
$P_i(t_1,\dots,t_i)$ with rational coefficients, such that:
\[
\chi_M(t)=t^n+\sum_{i=1}^nP_i(tr(M),\dots,tr(M^i))t^{n-i}.
\]
The polynomials
$P_i(t_1,\dots,t_i)$ are the ones which express the elementary
symmetric functions
$e_i(x_1,\dots,x_n)$ defined by
$\prod(t-x_i)=t^n+\sum_{i=1}^ne_i(x_1,\dots,x_n)t^{n-i}$ in terms of the Newton
functions $\psi_k:=\sum_ix_i^k$, i.e.
$e_i(x_1,\dots,x_n)=P_i(\psi_1,\dots,\psi_i).$\smallskip

At this point we can formally define, in an algebra with trace $R$, for every
element $a$ a formal $n-$characteristic polynomial:
\[
\chi_a^n(t):=t^n+\sum_{i=1}^nP_i(t (a),\dots,t (a^i))t^{n-i}.
\]
With this definition we obviously see that, given any element $a$, the element
$\chi_a^n(a)$ vanishes in every $n-$dimensional
representation or equivalently it is in the kernel of the universal map.
Thus we are led to make the following.

\begin{definition}
An algebra with trace $R$ is said to be an $n-$Cayley Hamilton algebra,
or to satisfy the $n^{th}$ Cayley Hamilton identity if:

1) $t(1)=n$.\footnote{this restriction could be to some extent dropped}

2) $\chi_a^n(a)=0, \ \forall a\in R$.
\end{definition}

It is clear that  $n-$Cayley Hamilton algebras form a category.  This category
has obviously free algebras. By definition the free $n$-Cayley-Hamilton algebra
$ F_n\langle x_i\rangle$
is the algebra generated freely by  $x_i$ and by traces of monomials
    modulo the trace ideal generated by evaluating in all possible ways
the $n^{th}$ Cayley
Hamilton identity. It is thus a quotient of the free algebra
$ F\langle x_i\rangle$ with trace.
Some remarks are in order, they are all standard from the theory of 
identities \cite{P},\cite{P3}.

One can {\it polarize} the Cayley-Hamilton identity getting a
multilinear identity  $CH(x_1,\dots,x_n)$, since we are in
characteristic 0 this identity is equivalent to the 1-variable
Cayley-Hamilton identity, $CH(x_1,\dots,x_n)$ has a nice
combinatorial description as follows.  given a permutation $\sigma$
of $n+1$ elements written into cycles as
$\sigma=(i_1,i_2,\dots,i_a)(j_1,j_2,\dots,j_b)\dots
(v_1,v_2,\dots,v_r)(u_1,u_2,\dots,u_s,n+1)$ set
$$\phi_\sigma(x_1,\dots,x_n):= t( x_{i_1}x_{i_2}\dots
x_{i_a})t(x_{j_1}x_{j_2}\dots x_{j_b})\dots t(x_{v_1}x_{v_2}\dots
x_{v_r})x_{u_1}x_{u_2}\dots x_{u_s}
$$
then
$$CH(x_1,\dots,x_n)=(-1)^{n}\sum_{\sigma\in S_{n+1}}
sgn(\sigma)\phi_\sigma(x_1,\dots,x_n).$$

{\bf Base change}  If $R$ is a $A-$algebra satisfying the
$n$-Cayley-Hamilton identity and $B$ is a commutative $A$
algebra,  then
$R\otimes_AB$ acquires naturally a $B-$linear trace for which it is
also an $n$-Cayley-Hamilton algebra.\smallskip

    By
construction the free trace algebra satisfying the
$n^{th}$ Cayley Hamilton identity in variables $x_i,\ i\in I$ has as
universal map to the algebra $C_n(\xi_1,\dots,\xi_m)$
of $GL(n)$-equivariant maps of matrices
$M_n(k)^I\to M_n(k)$.

The first main theorem of the theory is the
following:
\begin{theorem}\label{univ-mat} \begin{enumerate}
\item The universal map
\[
F_n\langle x_i\rangle\stackrel{i}\longrightarrow C_n(\xi_1,\dots,\xi_m)
\]
from the free  $n^{th}$ Cayley Hamilton algebra in variables $x_i$ to
the ring of
equivariant maps on
$m-$tuples of matrices is an isomorphism.

\item The trace algebra $T_n(\xi_1,\dots,\xi_m)$ of $C_n(\xi_1,\dots,\xi_m)$
is the algebra of invariants of $m-$tuples of matrices. As soon
as
$m>1$, it is the center of $C_n(\xi_1,\dots,\xi_m)$.

\item $C_n(\xi_1,\dots,\xi_m)$ is a finite $T_n(\xi_1,\dots,\xi_m)$ module.
\end{enumerate}
\end{theorem}
\proof
We sketch the proof, (cf. \cite{P3}, \cite{P4}, \cite{R}). We take an
infinite set of
variables $x_i,\
i=1,\dots,\infty$ we let the linear group
$G:=GL(\infty,k)$ act by linear transformations on variables on both sides.
The map $i$ is clearly $G$ equivariant. By standard
representation theory, in order to prove that $i$ is an isomorphism
it is enough
to check it on multilinear elements.  Now the space of
multilinear equivariant maps of matrices  in $m-$variables
$f:M_n(k)^m\to M_n(k)$
can be identified with the space of multilinear
invariant functions of matrices  in $m+1-$variables  $g:M_n(k)^{m+1}\to k
$ by the formula
$tr(f((\xi_1,\dots,\xi_m)\xi_{m+1})$. This last space
can be identified with the centralizer of $Gl(n,k)$ acting on the
${m+1}^{th}$ tensor
power of $k^n$. Finally this space is identified with
the group algebra of the symmetric group  $S_{m+1}$ modulo the ideal generated
by the antisymmetrizer on $n+1$ elements. Finally one has
to identify the element  of the symmetric group decomposed into cycles
$(i_1,\dots,i_k)\dots (s_1,\dots,s_l) (j_1,\dots,j_h,m+1)$  with
covariant map
$tr(\xi_{i_1}\dots \xi_{i_k})\dots tr(\xi_{s_1}\dots \xi_{s_l})\xi_{j_1}\dots
\xi_{j_h}$ and finally identify the antisymmetrizer with the
Cayley Hamilton identity and the elements in the ideal of the symmetric group
with the elements deduced from this identity in the free
algebra.  ii) follows easily from the previous description. For iii) one has
the estimate of Razmyslov that  $C_n(\xi_1,\dots,\xi_m)$ is
generated as
$T_n(\xi_1,\dots,\xi_m)$ module by the monomials in the $\xi_i$ of degree
$\leq n^2$. Conjecturally the right estimate is rather $\leq
{n+1 \choose 2}$.
\qed 

\noindent From this theorem, using a suitable method of Reynolds operators,
one gets (\cite{P5}):

\begin{theorem} If $R$ is a $n^{th}$ Cayley Hamilton algebra, the universal
trace preserving  maps $i_R,\overline i_R$
\[\begin{aligned}
    R&\stackrel{i_R}\longrightarrow  M_n(B_n(R))^{GL(n,k)}\\
{t}\downarrow &\quad\quad\quad    {tr}
\downarrow& \\
    T&\stackrel{\overline i_R}\longrightarrow   B_n(R)^{GL(n,k)}
\end{aligned}\]
     are   isomorphisms.
\end{theorem}

\section{Semisimple representations}

In this section and the next ones we will work on  finitely generated
algebras over an algebraically closed field $k$.  If
$A$ is a commutative  finitely generated algebra  over   $k$, we set
$V(A)$ to be the associated algebraic variety,  which
can be identified either to the maximal spectrum of $A$ or to the
homomorphisms $\phi:A\to k$.\smallskip

    Let us
look at $ F_n\langle x_i\rangle= C_n(\xi_1,\dots,\xi_m)$ in the case
that $k$ is
algebraically closed.  In this case we can apply geometric invariant theory.
First we analyze the case of the free algebra in
$m-$variables.

Here we have seen that the  coordinate ring of the $n-$dimensional
representations
of $ F\langle x_1 ,x_ 2,\dots ,x_m\rangle$ or
equivalently of  $ F_n\langle x_1 ,x_ 2,\dots x_m\rangle$ is the
coordinate ring of
the space of $m-$tuples of $n\times n$ matrices
$M_n(k)^m$.  The action of the linear group is by simultaneous
conjugation.  The
ring of invariants, which we have denoted by
     $T_n(\xi_1,\dots,\xi_m)$ is the coordinate ring of the {\it
quotient variety}
\[
V_{n,m}:=M_n(k)^m//Gl(n,k)
\]
the quotient map $M_n(k)^m\stackrel{\pi}\longrightarrow
V_{n,m}=M_n(k)^m//Gl(n,k) $
is surjective and each fiber
contains exactly one closed orbit.  By the analysis of
M. Artin (cf. \cite{A},\cite{P},\cite{P2}) we have that an $n-$tuple
$(A_1,\dots,A_m)$ of
matrices is in a
closed orbit if and only if it is {\it semisimple} in the sense that
the subalgebra
$k[A_1,\dots,A_m]\subset M_n(k)$ is semisimple.  Given any  $n-$tuple
$(A_1,\dots,A_m)$ of matrices the unique closed orbit contained in the
closure of its orbit is
constructed by taking a composition series of $k^n$
thought as a  $k[A_1,\dots,A_m]$ module and constructing the
associated graded
semisimple representation (whose isomorphism class is
uniquely determined).

It is interesting to analyze more closely  this picture. Let us use again
the notation $A_{n,m}:=k[M_n(k)^m]$ the coordinate ring of
the space of $m-$tuples of matrices. Given a point
$p\in V_{n,m}$ this is given by a maximal ideal $m_p$ of
$T_n(\xi_1,\dots,\xi_m)$, by the previous theory it corresponds to an
equivalence class of semisimple representations of  $
F\langle x_1 ,x_ 2,\dots ,x_m\rangle$, the closed orbit in the fiber
$\pi^{-1}(p)$.
An explicit   representation $\phi$ of $
F\langle x_1 ,x_ 2,\dots ,x_m\rangle$ in the fiber of  $p$ is given by a
maximal ideal  $M_\phi$ in the coordinate ring of matrices
lying over $m_p$ and the representation is given by the evaluation:

\[
\phi: F_n\langle x_i\rangle\stackrel{i}\longrightarrow
M_n[A_{n,m}]\longrightarrow
M_n[A_{n,m}/M_\phi]=M_n[k]
\]

Take now any finitely generated Cayley Hamilton algebra $R$,  and let
$T$ be its trace algebra.  $R= F_n\langle x_1
,x_ 2,\dots ,x_m \rangle/I$ is the
quotient of the free   Cayley Hamilton algebra  $ F_n\langle x_1
,x_ 2,\dots ,x_m \rangle$  modulo a trace ideal  $I$ and
correspondingly   $T$ is the
quotient of $T_n(\xi_1,\dots,\xi_m)$ modulo $I\cap T_n(\xi_1,\dots,\xi_m)$.

By Theorem 2.6 and functoriality we have a
commutative diagram:
\[
\begin{array}{lcll}
F_n\langle x_i\rangle &\stackrel{i\ \cong}\longrightarrow
&M_n[A_{n,m}]^{GL(n,k)}\longrightarrow &M_n[A_{n,m}]\\
\downarrow & & \downarrow& \downarrow\\
R&\stackrel{i_R\ \cong}\longrightarrow
&M_n[B_n(R)]^{GL(n,k)}\longrightarrow &M_n[B_n(R)]
\end{array}
\]

The ring $B_n(R)=A_{n,m}/I'$ need not be reduced, nevertheless it
defines a $GL(n,k)$ stable
subvariety of $M_n(k)^m$ made of the (trace) representations of
$R$, furthermore   $T= B_n(R) ^{GL(n,k)}$. Again explicitly a
homomorphism $\psi:A_{n,m}\to B_n(R)\to k$ gives
maximal ideals $M_\psi,M'_\psi$ and the commutative diagram:
\[
\begin{array}{lclcl}
F_n\langle x_i\rangle &\stackrel{i}\longrightarrow
&M_n[A_{n,m}]&\longrightarrow &
     M_n[A_{n,m}/M'_\psi]=M_n[k]\\
\downarrow & & \downarrow & & 1 \downarrow \\
R & \stackrel{i_R}\longrightarrow& M_n[B_n(R)]&\longrightarrow &
M_n[B_n(R)]/M_\psi]=M_n[k]
\end{array}
\]
If $\psi:A_{n,m}\to B_n(R)\to k$  gives a point in a closed orbit
then the corresponding representation of $R$ is
semisimple. It follows:
\begin{theorem}\label{Ct}  The algebraic variety associated to the
ring $T$ parametrizes isomorphism classes of (trace
compatible) semisimple representations of $R$.

\end{theorem}
Consider again $ \phi: C_n(\xi_1,\dots,\xi_m)= F_n\langle
x_i\rangle\stackrel{i}
\longrightarrow M_n[A_{n,m}]\longrightarrow
M_n[A_{n,m}/M_\phi]=M_n[k]$ which induces a map $\phi:
T_n(\xi_1,\dots,\xi_m)
\longrightarrow k$ with kernel some maximal ideal $m_p$.
The representation  $\phi$ factors   through
$$C_n(p):=C_n(\xi_1,\dots,\xi_m)/m_p
C_n(\xi_1,\dots,\xi_m).$$

     The algebra
$C_n(p)$ is by construction a finitely generated (over the trace
algebra) Cayley Hamilton algebra with trace and the trace takes
values in
$k$, so we need to start with analyzing this picture.

\begin{proposition}\label{CH} Let $R$ be a $n^{th}-$Cayley Hamilton algebra
with trace  values in
$k$ and finitely generated  over $k$, denote by $t$ the trace,
then\footnote{one could drop the hypothesis that $R$ is finitely generated
and still get 2),3) but not 1)}:
\begin{enumerate}
\item $R$ is finite dimensional over $k$.

\item The trace is 0 on the Jacobson radical $J$ of $R$. So $t$
factors through $R/J$ which is also
    a $n^{th}-$Cayley Hamilton algebra.

\item $R/J=\oplus_{i=1}^sM_{k_i}(k)$ and there exists positive integers $h_i$
with $n=\sum_{i=1}^sh_ik_i$ such that given $a_i\in M_{k_i}(k)$ and
$tr(a_i)$ the ordinary
trace, we have:
\[
t(a_1,a_2,\dots,a_s)=\sum_{i=1}^s h_itr(a_i)
\]
\end{enumerate}
\end{proposition}
\proof  The finite dimensionality follows from Theorem \ref{univ-mat} 3).
By the previous
discussion, if we present $R$ as a quotient of
$C_n(\xi_1,\dots,\xi_m)= F_n\langle x_i\rangle$ we have a commutative diagram:
\[
\begin{array}{lcl}
C_n(\xi_1,\dots,\xi_m) & \stackrel{tr}\longrightarrow &
T_n(\xi_1,\dots,\xi_m)\\
\downarrow & & p \downarrow\\
R & \stackrel{t}\longrightarrow & k
\end{array}
\]
thus a point $p$ in $V_{n,m}$ which  can be lifted  to a semisimple
representation of $C_n(\xi_1,\dots,\xi_m)$ factoring through $R$.

It is enough to see that the kernel $ Ker\phi$  of this semisimple
representation
$\phi:R\to M_n(k)$ is exactly the radical $J$ since then
one can use the ordinary theory of semisimple algebras and deduce that
$R/J=\oplus_{i=1}^sM_{k_i}(k)$ and an $n-$dimensional
representation of this algebra is of the form $\oplus h_ik^{k_i}$ with trace
exactly $t(a_1,a_2,\dots,a_s)=\sum_{i=1}^s h_itr(a_i) $.

Clearly
$J\subset Ker\phi$ so it suffices to show that $Ker\phi$ is a nilpotent ideal.

By construction we have that, for every element $a\in
Ker\phi$ the (formal) trace $t(a)$ is 0.  Therefore the Cayley
Hamilton polynomial
is $\chi_{ a}^n(t)=t^n$ hence every element $a\in
Ker\phi$ is nilpotent and thus the Kernel is nilpotent.
\qed

The previous proposition applies thus to the algebras
\[
C_n(p):=C_n(\xi_1,\dots,\xi_m)/m_p
C_n(\xi_1,\dots,\xi_m).
\]
We need to understand a special case, when $\phi$ is irreducible.

\begin{proposition} If $\phi$ is irreducible then
\[
C_n(p):=C_n(\xi_1,\dots,\xi_m)/m_p
C_n(\xi_1,\dots,\xi_m)=M_n(k).
\]
\end{proposition}
\proof
     By the previous proposition we have that $C_n(p)/J=M_n(k)$ where $J$ is the
radical, we have to show that $J=0$. If $u_i,\
i=1,\dots,n^2$ are elements of $C_n(\xi_1,\dots,\xi_m)$ whose images in
$C_n(p)/J=M_n(k)$ form a basis, we have that the matrix
$tr(u_iu_j)$ is invertible in the local ring $T_n(\xi_1,\dots,\xi_m)_{m_p}$
since the trace form on $M_n(k)$ is non degenerate.
Moreover these elements form also a basis for the ring of matrices
$M_n(A_{n,m})$
localized at $m_p$. Thus given any
$u\in C_n(p)$ we have
$u=\sum_ix_iu_i,\ x_i\in  [A_{n,m}]_{m_p}$. The $x_i$ can be computed by the
equations $tr(uu_j)=\sum_ix_itr(u_iu_j)$ which shows that
$x_i\in T_n(\xi_1,\dots,\xi_m)_{m_p}$ specializing modulo $m_p$ we see that the
classes of the $u_i$ are a $k-$basis of $C_n(p)$ hence
$C_n(p)=C_n(p)/J=M_n(k)$.
\qed

    From the previous analysis it follows that the set of points $p$ in which
$C_n(p)= M_n(k)$ is the open set defined by the non vanishing of at
least one of the {\it discriminants} $\det(tr(u_iu_j))$ where the elements
$u_i$ vary on all $n^2-$tuples of elements of
$C_n(\xi_1,\dots,\xi_m).$

This set is empty if and only if $m=1$. For $m=2$ one takes a
diagonal matrix with
distinct entries and the matrix of a cyclic permutation which
generate $M_n(k)$.

A more careful analysis of the previous argument shows that the
localized algebra
$C_n(\xi_1,\dots,\xi_m)_{m_p}$ is an Azumaya algebra over
$T_n(\xi_1,\dots,\xi_m)_{m_p}$. In the geometric language, the points of the
spectrum of  $T_n(\xi_1,\dots,\xi_m)_{m_p}$ where
$C_n(\xi_1,\dots,\xi_m)_{m_p}$  is an Azumaya algebra are exactly the points
over which the quotient map
$$M_n(k)^m\stackrel{\pi}\longrightarrow M_n(k)^m//Gl(n,k)$$
is a {\it principal bundle} over the projective linear group.
     In fact this can be viewed as a special case of M. Artin 
characterization of
Azumaya algebras by polynomial identities (cf. \cite{A},
\cite{P},\cite{S}).

     One easily obtains the following
\begin{corollary} Let $R$ be an $n-$Cayley Hamilton algebra with
trace  values in
$k$ and  Jacobson radical $J$.  If  $R/J= M_{n}(k)$ then $J=0$.
\end{corollary}
\vskip20pt
\subsection{Some categorical constructions.}
Let $R$ be a  $ n^{th}-$Cayley-Hamilton algebra finitely
generated algebra over an algebraically closed field $k$,  with trace 
$t$ and $t(R)=A$. $A$ is also a
finitely
generated algebra over   $k$. 
By Theorem  \ref{Ct} the (closed) points of $V(A)$
parametrize
semisimple representations of
dimension $n$ of $R$. Fix a positive integer $r$ and   change trace 
taking the new trace
$\tau=rt_{R/A}$.

\begin{proposition} \label{pm} Let $R$ be a  $ n^{th}-$Cayley-Hamilton algebra
with trace $t$ and $t(R)=A$. $A$ a finitely
generated algebra over an algebraically closed field $k$ and $V(A) $
the    reduced
variety of $Spec(A)$.

    The algebra $R $ with
trace  $\tau:=rt $  is an
$(rn)^{th}-$Cayley-Hamilton algebra. Given
      a point $p\in V(A)$ it determines  an $n-$dimensional semisimple
representation $M_p^t$ compatible with the trace $t$
and also  an $rn-$dimensional semisimple representation $M_p^\tau$
compatible with the trace $\tau$  we have:
$$M_p^\tau= r M_p^t=(M_p^t)^{\oplus  r}$$
\end{proposition}
\proof Let us show first that, if $R$ is a  $ n
^{th}-$Cayley-Hamilton algebra with trace $t$  the algebra
$R $ with trace  $\tau:=rt $  is an
$(rn)^{th}-$Cayley-Hamilton algebra.  We know that $R$ embeds in a
trace compatible way into $n\times n$ matrices over a
commutative ring $A$, then the natural diagonal embedding of
$M_n(A)$ into $M_{rn}(A)$  gives the claim.

Now we have that $V(A)$ parametrizes also semisimple representations
(compatible with the new trace) of dimension $rn$. We have an
obvious map from the variety of $n$ dimensional representations compatible
with the reduced trace to the variety of $rn$ dimensional
representations compatible with   trace $\tau$
it is simply the map that associates
to a representation $M$ its direct sum $M^{\oplus r}$.

    From this the statement is clear.\qed

Let now $R_1$ and $R_2$  two trace algebras over $A$  which are
Cayley-Hamilton for two integers $n_1,n_2$, and
$t(R_1)=t(R_2)=A$,  then, a point $p\in V(A)$ determines a semisimple
representation $M^1_p$ of dimension $n_1$ of $R_1$
and a semisimple representation $M^2_p$ of dimension $n_2$ of $R_2$:
\begin{proposition}\label{ds}  The algebra $R:=R_1\oplus R_2$ with
trace $t(r_1,r_2):= t(r_1)+t(r_2)$ is an
$n^{th}-$Cayley-Hamilton  for   $n:=n_1+n_2$, and
$t(R)=A$,   a point $p\in V(A)$ determines also a semisimple
representation $M_p$ of dimension $n=n_1+n_2$ of $R$ and:
$$M_p=M_p^1\oplus M_p^2 $$
\end{proposition}

The proof is similar to that of the previous proposition and it is
omitted.\vskip20pt

Finally let now $R_1$ and $R_2$  two trace algebras over $k$  which
are Cayley-Hamilton for two integers $n_1,n_2$, and
$t(R_1)=A_1,\ t(R_2)=A_2$,  then, a point $p  \in V(A_1)$ determines
a semisimple representation $M^1_p$ of dimension $n_1$ of $R_1$
and a   point $q  \in V(A_2)$ determines a semisimple representation
$M^2_q$ of dimension $n_2$ of $R_2$:
\begin{proposition}\label{tp}  The algebra $R:=R_1\otimes R_2$ with
trace $t(r_1\otimes r_2):= t(r_1)\otimes t(r_2)$ is an
$n^{th}-$Cayley-Hamilton  for   $n:=n_1 n_2$, and
$t(R)=A:=A_1\otimes A_2$,   a point $(p,q)\in V(A)=V(A_1)\times
V(A_2)$ determines also a semisimple representation $M_p$ of
dimension
$n=n_1 n_2$ of
$R$ and:
$$M_p=M_p^1\otimes M_p^2 $$
\end{proposition}

Again the proof is similar.\vskip20pt
\section{The reduced trace}

Let us recall that a prime ring $R$ is a ring in which the product of two
non-zero ideals is non-zero.  Let $R$ be a prime algebra over a
commutative ring $A$ and assume that $A\subset R$ and $R$ is an $A-$module of
finite type. One easily sees that:

\begin{enumerate}
\item $A$ is an integral domain.

\item $R$ is a torsion free module.

\end{enumerate}

If $F$ is the field of fractions of $A$ then $R\subset R\otimes_AF$ and
$S:=R\otimes_AF$ is, by  a  Theorem of Wedderburn a (finite
dimensional) simple algebra isomorphic to $M_k(D)$ where $D$
is a finite dimensional division ring.

If $Z$ is the center of $S$, it is also the center of $D$ and
$\dim_ZD=h^2$; moreover,
if $\overline Z$ is
an algebraic closure of $Z$ we have
$M_k(D)\otimes_Z\overline Z=M_{hk}(\overline Z)$, if the finite extension
$Z\supset F$ is separable, as
happens in characteristic 0 and if $p:=[Z: F]=\dim_FZ$ we also have
\[
S\otimes_F\overline Z=M_k(D)\otimes_F\overline Z=M_{hk}(\overline Z)^{\oplus p}
\]

We define the number $hkp$ to be the {\it degree} of $S$ over $F$ and let:
\[
hkp:=[S:F],
\]

If $a\in S$ we have that $a\otimes 1\in M_k(D)\otimes_F\overline Z$
is a $p-$tuple
of matrices $(A_1,\dots, A_p)$ one defines the {\it reduced
trace} of $a$ to be the sum
\[
t(a)=t_{S/F}(a):=\sum_{i=1}^ptr(A_i),
\]
A standard argument of Galois theory shows that $t(a)\in F$,  in
characteristic 0
this can be even more easily seen as follows.

Consider the $F-$linear operator $a^L:S\to S,\ a^L(b):=ab$ let us
compute its trace.
This can be done in $M_k(D)\otimes_F\overline
Z=M_{hk}(\overline Z)^{\oplus p}$ where $a^L= (A_1^L,\dots, A_p^L)$ and so
\[
tr(a^L)=\sum_{i=1}^ptr(A_i^L)=hk\sum_{i=1}^ptr(A_i)\quad\implies
t_{S/F}(a)=\frac{1}{hk}tr(a^L)
\]

\begin{theorem}  If $S=R\otimes_A F$ as before and $A$ is integrally closed we
have that the reduced trace $t_{S/F}$ maps $R$ into $A$, so
we will denote by $t_{R/A}$ the induced trace.

The algebras $R, S$
with their reduced trace are
$n-$Cayley Hamilton algebras of degree
$n=hkp=[S:F]=[R:A]$, (we set $[R:A]:=[S:F]$).
\end{theorem}
\proof We have a natural representation of $M_{hk}(\overline
Z)^{\oplus p}$ by $hkp$ matrices for which the reduced trace
is the trace.  If  $R$ is a finite $A-$module, it is easy to see that
the reduced trace of an element
of $r$ is integral over $A$. If $A$ is integrally closed then the
trace takes values in $A$.\qed

The importance of the reduced trace comes from the next result. $R,A$
are as before,
$n:=[R:A]$ the degree and we are assuming characteristic 0:
\begin{theorem}\label{Rtrthm}  If $\tau:R\to A$ is any trace for which  $R $
is an
$m-$Cayley Hamilton algebra then there is a positive integer $r$ for which:
\[
m=rn,\qquad \tau=r\,t_{R/A}
\]
\end{theorem}
\proof  Clearly $\tau$ extends to a trace on $S$ with values in $F$
for which  $S $
is an $m-$Cayley Hamilton algebra.

Let $G$ be a finite Galois extension of $F$ for which
$S\otimes_FG=M_{hk}(G)^{\oplus p}$, $\tau$ extends to a $G-$valued trace
on $S\otimes_FG=M_{hk}(G)^{\oplus p}$ which is invariant under the Galois group
and  for which  $S\otimes_FG $
is an
$m-$Cayley Hamilton algebra.   Passing to the algebraic closure we can now
apply statement 3 in Proposition \ref{CH} where we know that, by
invariance under the Galois group, all the integers $h_i$ must be
equal to some
positive integer $r$. The formula follows from the definitions.
\qed

\section{The unramified locus and restriction maps}
\subsection{The unramified locus}

Let us go back to the previous setting. $R$ a prime algebra over $A$,
$F$ the field
of fractions of $A$, $S=R\otimes_AF$ and finally $Z$ the
center of $S$. Let now $B:=R\cap Z$ be the center of $R$.  If we
further assume that
$B$ is integrally closed we have the reduced trace
$t_{R/B}$ and the formulas:
\[
t_{R/A}=t_{B/A}\circ t_{R/B},\quad [R:A]=[R:B][B:A]
\]

Let us assume now that $A$ is a finitely generated algebra over an
algebraically
closed field $k$ and $V(A) $   the   associated affine
variety    parame\-trizing semisimple representations
(compatible with the reduced trace $t_{R/A}$) of $R$ of dimension
$n=[R:A].$ Since we are assuming that $R$ is a finite $A$ module it
follows that also
$B$ is a finite $A$ module. Then $B$ is a finitely
generated algebra over   $k$ and its   associated affine variety $V(B) $
     parametrizes semisimple
representations  (compatible with the reduced trace $t_{R/B}$) of $R$ of
dimension $m=[R:B].$  Moreover we can also use $V(A)$ to parametrize
semisimple representations
(compatible with the reduced trace $t_{B/A}$) of $B$ of
dimension $p=[B:A].$  Finally the inclusion $A\subset B$
defines a morphism of algebraic varieties $\pi:V(B)\to V(A)$ of degree $p$.
We want to
put all these things together.

Given a point $Q\in V(A)$ denote by $N_Q$ the corresponding $mp$ dimensional
semisimple representation
of $R$. Given a point $P\in V(B)$ denote by $M_P$
the corresponding $m $ dimensional semisimple representation of $R$
.\smallskip

First of all an irreducible
representation of $B$ is 1-dimensional and corresponds to a point
$P\in V(B)$, a
semisimple representation corresponds to a positive cycle $\sum
h_iP_i$ of degree $p=\sum_ih_i$. Proposition \ref{CH} implies:

\begin{proposition} Given a point $Q\in V(A)$ we have for the
associated semisimple
representation
$\sum_{i=1}^s h_iP_i$ of $B$, that the points $P_i$ are exactly the
points in the
fiber
$\pi^{-1}(Q)$. So we may identify formally $\sum_{i=1}^s h_iP_i$ with
the cycle  $[\pi^{-1}(P)].$
\end{proposition}

In general a fiber need not have exactly $p$ points but it can have
$s\leq p$ points.

In terms of algebras, $Q$ corresponds to a maximal ideal $m$ of $A$
and the points
$P_i$ to the maximal ideals of $B/mB$. This is a finite
dimensional commutative algebra and so  $\overline B:=B/mB=\oplus_{i=1}^s B_i$
where $B_i$ is a local ring supported in the point $P_i$. Let
$\overline n_i$ be the maximal ideal of $B_i$ and $n_i$ the
corresponding maximal ideal of $B$. We have $B_i/\overline
n_i=k$ and again Proposition
\ref{CH}
implies that the trace $t_{B/A}$ induces a trace $t_{\overline B/k}$
on $\overline B$ decomposes
as the sum of {\it local factors} $h_it_i$ where $t_i:B_i\to
B_i/\overline n_i=k$ is the
projection.  In other words, if $e_i$ is the idempotent, unit of
$B_i$  we have $t_{\overline B/k}(e_i)=h_i$.

\begin{remark}  If $B$ is a projective $A-$module of rank  $n$ the
reduced trace  $t_{B/A}(b)$ is just the trace of the
linear map $x\to bx$. In this case  $\dim_k\overline B =n,\quad\dim_k B_i=h_i$

\end{remark}

Passing to the algebra $R$ we have  a direct sum decomposition.
\[
\overline R:=R/mR=R\otimes_AA/m=
R\otimes_B(B\otimes_AA/m)=\oplus_{i=1}^sR\otimes_B
B_i=\oplus_{i=1}^s\overline R_i
\]

We pass to the traces: we know that $t_{R/A}=t_{B/A}\circ t_{R/B}$,
and modulo $m$
we get traces  $t_{\overline R/k}=t_{\overline
B/k}\circ t_{\overline R/\overline B}$. If $e_i\in B_i$ is the
idempotent identity of $B_i$ we have $\overline
R_i=\overline Re_i$  thus   $ t_{\overline R/\overline B}$ restricts
to $\overline R_i$ to a $B_i$ trace
     and  $t_{\overline R/\overline
B}=\oplus_{i=1}^st_{\overline R_i/B_i}$.   Similarly $t_{\overline
B/k}=\oplus_{i=1}^st_{ B_i/k}$ and
$$t_{\overline R/k}= \oplus_{i=1}^s t_{  B_i/k}\circ t_{\overline R_i/B_i}$$

The algebra $R_i$ with the trace $t_{\overline R_i/B_i}$  is an
$m=[R:B]$-Cayley Hamilton algebra,
     if $\overline n_i$ is the maxiamal ideal of $B_i$ the unique point
of $Spec(B_i)$ given by $B_i\to B_i/\overline n_i=k$
corresponds to some semisimple representation  $M_i$, as $R$
representation it is $M_{P_i}$ where $P_i$ corresponds to
the maximal ideal $n_i$, since we have that
$R/n_iR=\overline R_i/\overline n_i\overline R_i$. If we
denote by
$\overline t_{\overline R_i/B_i}$ the image of
$t_{\overline R_i/B_i}$ modulo  $\overline n_i$ we have
$$h_i\overline t_{\overline R_i/B_i}=t_{\overline B_i/k}\circ
t_{\overline R_i/B_i}
    $$ it follows that the semisimple representation of  $\overline R_i$
relative to the trace $t_{\overline B_i/k}\circ
t_{\overline R_i/B_i}$  is $h_iM_i$.
Thus the semisimple representation of $\overline R$ relative to the
trace $t_{\overline R/k}= \oplus_{i=1}^s t_{
B_i/k}\circ t_{\overline R_i/B_i}$ is $\oplus_ih_iM_{P_i}$, (cf.
\ref{pm}, \ref{ds}).\smallskip

We have proved:
\begin{theorem}
Given a point $Q\in V(A)$  and its cycle $\sum_ih_iP_i$ in $V(B)$ we have
\[
N_Q=\oplus_i h_iM_{P_i}
\]
\end{theorem}
Let $R$ be an algebra with trace
finitely generated over an algebraically closed field $k$
satisfying the $n^{th}$ Cayley Hamilton identity with the trace algebra
$t(R)$. If $P\in V(t(R))$ with maximal ideal
$m_P$ we ask when $R(P):=R\otimes_{t(R)}k=R/m_PR$ is a semisimple
algebra?

The answer is implicit in Proposition \ref{CH}. The trace map
for $R$ induces the
trace $t:R(P)\to k$ for $R(P)$ and the bilinear trace form $t(ab)$.
It follows immediately from \ref{CH} that:
\begin{proposition}
The radical $J$ of $R(P)$ is the kernel of the trace form $t(ab)$.
\end{proposition}

Let us see what is the meaning of this statement
in the case $R\supset A$ is a prime algebra over $A$
finitely generated algebra over an algebraically closed field
$k$, $F$ the field of fractions of $A$, $S=R\otimes_AF$,   $Z$ the
center of $S$,
$B:=R\cap Z$ be the center of $R$ with    $B$   integrally
closed.  If $m:=[R:B]$ the map ${t_{R/B}\over m}$ is a projection on
$B$ so that $R=B\oplus R^0$ and $B$ is a direct
summand.

With the previous notations
      $V(A) $   the   associated affine
variety     parametrizes semisimple representations
(compatible with the reduced trace $t_{R/A}$) of $R$ of dimension
$n=[R:A].$   $B$ is a finitely
generated algebra over   $k$ and its   associated affine variety $V(B) $
      parametrizes semisimple
representations  (compatible with the reduced trace $t_{R/B}$) of $R$ of
dimension $m=[R:B].$  Denote by $\pi:V(B)\to V(A)$.\smallskip

If $Q\in V(A)$ corresponds to a maximal ideal $m_Q$ and the algebra
$R(Q):=R\otimes_AA/m_Q$ is
semisimple we have that $R\otimes_AA/m_Q=B\otimes_AA/m_Q\oplus
R^0\otimes_AA/m_Q$, since $B(Q):=B\otimes_AA/m_Q$ is in the
center of $R(Q):=R\otimes_AA/m_Q$ we have that
$B(Q)$ is semisimple, in other words the scheme theoretic fiber
$\pi^{-1}(Q)$ is reduced, $B\otimes_AA/m_Q=\oplus_iB/n_i$
and $B/n_i=k$ is a   point
$P_i$ in the fiber of
$Q$. We also have $R(P_i)=R\oplus B/n_i$ and
$R\otimes_AA/m_Q=\oplus_i R(P_i)$ hence $R(P_i)$ is semisimple. The
converse
is also clear.\smallskip

    The commutative algebra $B(Q)$ is semisimple if and only if it is
reduced, i.e. the fiber of $Q$
under the map $\pi$ is reduced, which in our case implies that $\pi$
is {\'e}tale in
the points of this fiber. Now we  know that
$R$ is a finite module over
$B$ and its generic dimension is $m^2$. The dimension of
$R(P)$ over $k=B(P),\ P\in V(B)$ is a semicontinuous function and we
always have
$\dim_{B(P)}R(P)\geq m^2$. If $R(P)$ is not simple of dimension
$m^2$ from \ref{CH}, 3) follows that, if $J$ is the radical of $R(P)$ we have
$\dim_kR(P)< m^2$  hence we have a dichotomy, either $R(P)=M_m(k)$
that is to say that $P$ corresponds to an irreducible representation,
or $R(P)$ is
not semisimple.
\begin{proposition}
Let $W^0$ be the open set of $V(B)$ made of points $P$ where $R(P)=M_m(k)$
(the irreducible representations), let $V^0$ be the maximal open set of
$V(A)$ with $\pi^{-1}(V^0)\subset W^0$ and $V^1$ be the  open set of
$V(A)$ where
$B(Q)$ is reduced, then the set of points of $V$ where $R(Q)$ is semisimple is
$V^0\cap V^1$.
\end{proposition}
The set of points $V^0\cap V^1$ of $V$ where $R(Q)$ is semisimple is called the
{\it unramified locus} of the $A$ algebra $R$.

    From our analysis it follows that, if $m=[R:B],\ p=[B:A]$ over a 
point $Q$ in
the unramified locus
\[
N_Q=\oplus_{i=1}^pM_{P_i}
\]
decomposes as the direct sum of the $p$
irreducible representations $M_{P_i}$ supported
at the
$p$ distinct points of the fiber
$\pi^{-1}(Q)$.

\subsection{Restriction maps}\label{restmaps}

We come now to the final application of the previous theory. The
setting we have
in mind appears naturally for quantum groups at roots of 1 and
their subgroups.\smallskip

We need a first Lemma. Given a prime algebra $R$ finite over $A$ with
center $Z$, let $F$ be the quotient field of $A$
and $F\subset G$ an extension field.
\begin{lemma} The following are equivalent:

i)  $R\otimes_AG$ is a simple algebra.

ii) The algebra $Z\otimes_AG$ is a field.
\end{lemma}
\proof Let $F$ be the quotient field of $A$ we have that
$S:=R\otimes_AF$ is a simple
algebra with center the field
$W:=Z\otimes_AF$ and that
$R\otimes_AG=(R\otimes_AF)\otimes_FG=S\otimes_W(W\otimes_FG) $. Since
$S$ is a simple algebra with center $W$         it is well
known and easy that
$S\otimes_W(W\otimes_FG)
$ is simple if and only if $ W\otimes_FG $  is a field.  Finally  $ W\otimes_FG
=(Z\otimes_AF)\otimes_FG=Z\otimes_A G.$
\qed

Assume that we have two prime algebras $R_1\subset R_2$ over two
commutative rings
$A_1\subset A_2\subset R_2$. Assume as in the previous
paragraph that each $R_i$ is finitely generated as $A_i$ module and
that the two rings $A_i$ are integrally closed. We thus have the two
reduced traces $t_{R_i/A_i}$, we want to discuss the compatibility of
these traces.
In general one can see by simple examples that there is no
compatibility. Let us thus make the basic assumption of {\it 
compatibility} (with trace).
We let $F_i$ be the quotient field of $A_i$  and consider
$S_i:=R_i\otimes_{A_i}F_i$.
\begin{lemma} Given two prime algebras
$R_1\subset R_2$ over the rings
$A_1\subset A_2\subset R_2$  with
$Z_1$ the center of
$R_1$. Assume $R_1$ is finite over $A_1$. The following two 
conditions are equivalent:

   i)  $R_1\otimes_{A_1}F_2$ is a simple algebra.

ii) The algebra
$Z_1\otimes_{A_1}F_2$ is a field.

In this case  the map $i:
R_1\otimes_{A_1}F_2=S_1\otimes_{F_1}F_2\to S_2$ is injective.
\smallskip

These conditions are satisfied if:

iii) The algebra
$Z_1\otimes_{A_1}A_2$ is a domain.

\end{lemma}
\proof The equivalence of the first two conditions is the content of 
the previous Lemma. It is clear
that iii) implies ii) since if $Z_1\otimes_{A_1}A_2$ is a domain, 
$Z_1\otimes_{A_1}F_2$ is its
quotient  field.

\begin{definition}\label{comp} We
say that the two algebras
$R_1\subset R_2$ are compatible with
$A_1\subset A_2\subset R_2$ if  the previous two equivalent 
conditions are satisfied.

\end{definition}
In the examples which we will study we will usually verify iii).

    \begin{remark} If $R_2$ is a domain then  $R_1\subset R_2$ is 
compatible with
$A_1\subset A_2\subset R_2$ if and only if  the map $i:
R_1\otimes_{A_1}F_2=S_1\otimes_{F_1}F_2\to S_2$ is injective.
\end{remark}
\proof
    If $R_2$ is a domain so is $S_2$ and so, if $i$ is injective we must
have that  $Z_1\otimes_{A_1}F_2$ is a
field.\qed

\begin{example} 1) If the extension $F_1\subset F_2$ is unirational
then $Z_1\otimes_{A_1}F_2$ is a field.

2) If $A_1=Z$ is the center of $R_1$ then  $Z_1\otimes_{A_1}F_2=F_2$
is a field.
\end{example}
    1) In fact we have
$F_2\subset F_1(t_1,\dots, t_m)$ so $Z_1\otimes_{A_1}F_2\subset
Z_1\otimes_{A_1}F_1(t_1,\dots,
t_m)=(Z\otimes_{A_1}F_1)(t_1,\dots, t_m)$ is a field.

\begin{theorem}\label{ress} Given two compatible algebras $R_1\subset
R_2$ we have that
for a positive integer $r$:
$$r [R_1:A_1]=  [R_2:A_2],\qquad r\,t_{R_1/A_1}=t_{R_2/A_2}
\quad\text{ on $R_1$}$$
\end{theorem}
\proof By the hypotheses made one can reduce the computation to the
two algebras $ S_1\otimes_{F_1}F_2\subset S_2$ over $F_2$. In this
case we know that the reduced trace  $t_{S_2/F_2} $ restricted to  $
S_1\otimes_{F_1}F_2$ makes it a Cayley Hamilton
trace algebra. Since by assumption $ S_1\otimes_{F_1}F_2$ is simple,
one can then apply Theorem 4.2.
\qed

Let us now assume to be in the geometric case in which $A_1,A_2$ are further
assumed to be finitely generated over an algebraically closed
field $k$. If $V(A_1),V(A_2)$ are the two associated affine varieties
parametrizing
semisimple representations we have an induced map $\pi:V(A_2)\to
V(A_1)$. If $Q\in V(A_2)$  and $M$ is a representation of $R_2$ over
$Q$ we then see
that $M$ is also a representation of $R_1$ over $\pi(Q)$ but
for $r$ times the reduced trace.  If $M$ is semisimple as $R_2$
module it may well be that it is not semisimple as $R_1$
module.

\begin{theorem}\label{rest} Given two compatible algebras $R_1\subset
R_2$ as before,
$Q\in V(A_2)$. $M_Q$ the corresponding semisimple representation
of $R_2$ of dimension  $ [R_2:A_2]$, $M_{\pi(Q)}$ the corresponding semisimple
representation
of $R_1$ of dimension  $ [R_1:A_1].$ We have that the restriction of
$M_Q$ to $R_1$
is a trace representation for $r\,t_{R_1/A_1}$, its
associated semisimple representation is $rM_{\pi(Q)}=M_{\pi(Q)}^{\oplus r}.$

If $ {\pi(Q)}$ lies in the unramified locus of $R_1$ (as $A_1$ algebra) we have
that the restriction of $M_Q$ to $R_1$ is the semisimple
representation
$rM_{\pi(Q)} .$
\end{theorem}
    \proof Everything follows from the previous discussions except the
last point.
Let $m\subset A_2$ be the maximal ideal associated to $Q$
and $m':=m\cap A_1$, by definition of the unramified locus    the algebra
$R_1/m'R_1 $ is a semisimple algebra for which every  representation
is semisimple.
\qed

    \subsection{Examples}

In this section we collect examples from quantum groups.

One class of examples is obtained by taking a quantum group $R$ at
roots of 1, where $A$ is a Hopf subalgebra coordinate
ring of an algebraic group $G$. In this case $R,\ R\otimes_{\mathbb
C}R$ are domains and we need to prove that:
\begin{theorem}\label{com} The
comultiplication  $\Delta:  R\to R\otimes_{\mathbb C}R$ is compatible
with $\Delta:  A\to A\otimes_{\mathbb C}A.$
\end{theorem}
\proof  By the
previous lemma it is enough to show that, setting $\mathcal Z$ the
center of $R$, we have that $\mathcal
Z\otimes_{\Delta(A)}(A\otimes_{\mathbb C} A)$ is a domain.\smallskip

Let us use some geometric language.  $A$ is the coordinate ring of a
connected algebraic group $G$ and the map
$\Delta:A\to A\otimes A$ is the comorphism associated to the
multiplication $G\times G
\stackrel{\mu}\longrightarrow  G$.  Let $\nu:G\times G\to G\times G$
be defined by $\nu(x,y):=(xy, y),$
    clearly $\nu$ is an isomorphism and $\mu$ can be identified to $\nu$
composed by the first projection.  Thus
$\nu^*:A\otimes A\to A\otimes A$ maps $\Delta(A)\to A\otimes 1$.
Using this isomorphism we see that:
$$\mathcal
Z\otimes_{\Delta(A)}(A\otimes_{\mathbb C} A)\cong (\mathcal
Z\otimes_{ A } A)\otimes_{\mathbb C} A=  \mathcal
Z \otimes_{\mathbb C} A$$
is a domain.\qed
\medskip

The second class of examples we have in mind is when $R_2$ is a Hopf
algebra and $R_1,A_1,A_2$ Hopf subalgebras.  In this
case $A_2$ is the coordinate ring of an algebraic group  $G_2$ and
$A_1$ that of a quotient  group.  In the case  where
$G_2$ is solvable, as for the quantized enveloping algebras, we will
always get that the extension $F_1\subset F_2$ is
rational and so we can also conclude.

For instance for $R_2=U_q(\mathfrak g)$  the quantized enveloping
algebra of a semisimple group and $R_1=U_q(\mathfrak b^+)$ we
have that $A_2$ is the coordinate ring of the {\it dual} group
$U^-\times T\times U^+$ and  $A_1$   the coordinate ring
of the quotient group $  T\times U^+$  hence the rationality statement.

\subsection{Cayley-Hamilton Hopf algebras}
We formalize the previous discussion as follows:
\begin{definition} A
Cayley-Hamilton Hopf algebra is a Hopf algebra such that:
\begin{itemize}
\item it is a Cayley-Hamilton algebra
\item the trace subalgebra is a Hopf subalgebra
\end{itemize}
\end{definition}

In the next section we will see plenty examples of such Hopf algebras
that are given by quantized universal enveloping algebras at roots of
unity.

Let $R$ be a Cayley-Hamilton Hopf algebra
with the trace subalgebra $A:=t(R)$. Assume $R$ is prime and a finite
$A$ module. Let $Z\supset A$ be the center of $R$, set
$m:=[R:A],\ n=[R:Z],\ p:=[Z:A]$  so that
$m=np$.  For a point $x\in V(A)$ (resp. $P\in V(Z)$) denote by $N_x$
(resp. $M_P$) the corresponding semisimple $m-$dimensional
representation  (resp. $n-$dimensional).

Assume that
$ A $ is finitely generated over an algebraically closed field $k$,
so that $V(A),V(Z)$ are affine algebraic varieties,  and let
$\pi:V(Z)\to V(A)$ be the corresponding map of varieties. The
comultiplication on $A$ defines an associative binary operation on $V(A)$.
The antipode defines the
inverse operation for this operation on $V(A)$, so $V(A)$ is an
algebraic group.

   From \ref{rest}, and \ref{com}  we see that:
\begin{proposition}\label{cgd1} If $x,y\in V(A)$ and $N_x,N_y$ are
the corresponding semisimple representations then
the semisimple representation associated to $N_x\otimes N_y$ is $m N_{xy}$.
\end{proposition}

    We will say that the pair of points $x,y\in V(A)$ is generic if
both points   and  their product in $V(A)$ lie in the unramified
locus. Such pairs of points form a Zariski open subvariety
in $V(A)\times V(A)$.

    For each point $Q\in V(Z)$ in the fiber of,
either $x,y,xy$ the corresponding representation $M_Q$ is irreducible.
$$N_x=\oplus_{P\in \pi^{-1}(x)}M_P,\  N_y=\oplus_{Q\in
\pi^{-1}(y)}M_Q,\  N_{xy}=\oplus_{R\in
\pi^{-1}(xy)}M_R.
$$  From \ref{cgd1} we get

\begin{equation}\label{cgdec}N_x\otimes N_y=mN_{xy},\quad
\bigoplus_{P\in \pi^{-1}(x),\ Q\in
\pi^{-1}(y)}\!\!\!\!\!\!\!\!\!M_P\otimes M_Q =\bigoplus_{R\in
\pi^{-1}(xy)}mM_R,
\end{equation}
    Let $V$ and $W$ be two irreducible representations of the CH-Hopf
algebra $R$ such
     that the restrictions of $V$ and $W$ to $A$ are given by
the action of
characters $x,y\in V(A)$ respectively. Thus $V=M_P,\ W=M_Q$ where
$P\in \pi^{-1}(x),\ Q\in \pi^{-1}(y)$.  The restriction to $A$ of the
tensor product
$V\otimes W$ has the same property,  with central
character $xy$ .

If the pair $x,y\in V(A)$ is generic the tensor product
$M_P\otimes M_Q$ is semisimple as an $R$-module and we have the:
\begin{theorem}\label{CG} Clebsch-Gordan decomposition (cf. equation
\ref{cgdec}):
\begin{equation}\label{cgdecomp}
M_P\otimes M_Q\simeq \oplus_{R\in \pi^{-1}(xy)} M_R^{\oplus
h_R^{P,Q}},\quad  \sum_Rh_R^{P,Q}=n,\ \sum_{P,Q}h_R^{P,Q}=m
\end{equation}\end{theorem}
For quantized enveloping algebras at roots of 1 we will prove the 
stronger statement that all the
multiplicities $h_R^{P,Q}$ are equal.
\vskip20pt

\section{Quantized universal enveloping algebras at roots of 1}
\subsection{The definition $\U_\e$}
Let ${\g}$ be a simple Lie algebra  of rank $n$ with the root system
$\Delta$. Denote by $Q,P$ the root and weight
lattice. Fix simple roots $\alpha_1,\dots ,\alpha_n\in\Delta_+$ and denote by
$(a_{ij})^r_{i,j=1}$ the corresponding Cartan matrix. Denote by
$d_i$   the length of the $i$-th simple
root.

For an odd positive integer $\ell$ denote by $\e$ a primitive
root of 1 of degree $\ell$ (in case of components of type $G_2$ we also need to
restrict to $\ell$ prime with 3).

For any lattice $Q\subset \Lambda\subset P$ we have a quantized
universal enveloping
algebra
$U_\e^\Lambda({\g})$. It  is the associative algebra with 1 over
$\mathbb C$ generated by
$K_\mu$, $\mu\in \Lambda$, and
$E_i,F_i$, $i=1,\dots ,n$ with defining relations:
\begin{eqnarray*}
K_\mu K_\nu & =& K_\nu K_\mu \ , K_\mu K_{-\mu}=1 \quad
K_0=1 \ , \\ K_\mu E_i &= & \e^{\alpha_i(\mu)}E_iK_\mu \\
K_\mu F_i &=& \e^{-\alpha_i (\mu)}F_iK_\mu  \ , \qquad
E_iF_j -F_jE_i = \delta_{ij}
(K_{\alpha_i}-K_{\alpha_i}^{-1})/(\e_i-\e_i^{-1}) \ , \\
&&\sum^{1-a_{ij}}_{k=0} (-1)^k
\left[ \begin{array}{c} 1-a_{ij}\\ k\end{array}\right]_{\e_i}\,
E_i^{1-k-a_{ij}} E_j\, E_k^k =0 \ , \quad i\neq j\\
&&\sum^{1-a_{ij}}_{k=0} (-1)^k
\left[ \begin{array}{c} 1-a_{ij}\\ k\end{array}\right]_{\e_i}\,
F_i^{1-k-a_{ij}} F_j\, F_k^k =0 \ , \quad i\neq j
\end{eqnarray*}
Here \ $\e_i=\e^{d_i}$,
\[
\left[ \begin{array}{c} m\\h\end{array}\right]_\e\,
= \frac{[m]_\e!}{[m-h]_\e![h]_\e!} \ ,\qquad
[h]_\e!=[h]_\e\dots [2]_\e[1]_\e \ , \qquad
\left[h\right]_\e = \frac{\e^h-\e^{-h}}{\e-\e^{-1}} \ .
\]
The map $\Delta$ acting on generators as
\begin{eqnarray}\label{comult}
\Delta K _\mu &=& K_\mu\otimes K_\mu \ , \\
\Delta E_i &=& E_i\otimes 1 + K_{\alpha_i}\otimes E_i \ ,\\
\Delta F_i &=& F_i\otimes K_{\alpha_i}^{-1} + 1\otimes F_i
\end{eqnarray}
extends to the homomorphism of algebras $\Delta:\U_\e\to
\U_\e\otimes \U_\e$. Here we used the notation $K_i=K_{\alpha_i}$.
The pair $(\U_\e,\Delta)$ is
a Hopf algebra with the counit $\eta(L_\mu)=1$,
$\eta(E_i)=\eta(F_i)=0$.

For $\Lambda=P$ we have the {\it simply connected} quantized algebra,
denoted by
$\U^s_\e$ or simply
$\U_\e$, for $\Lambda=Q$ we have the adjoint form denoted by $\U^a_\e$, the
definitions hold also if instead of $\epsilon$ we have a $q$ generic.

We will denote by $\U_\e^\pm$ the subalgebras of $\U_\e$
generated by $E_i$ and $F_i$ respectively. The subalgebra generated
by $K_i$ will be denoted $\U_\e^0$.

\subsection{PBW basis and the structure of the center}
One can introduce a monomial basis in
the algebras $\U_\e^\pm$ that is the analog of the
Poincare-Birkhoff-Witt basis. We will call it PBW basis.
In order to describe this basis we first should introduce
root elements $E_\alpha\in \U^+_\e$, $F_\alpha\in \U^-_\e$.
This can be done ( see \cite{L}
for details) by choosing a convex ordering on positive roots.

If $\beta(1)>\dots >\beta(N)$ is the convex ordering of positive
roots $\Delta_+$ ( here $N=|\Delta_+|$) then we choose PBW bases as
follows. For $\U^+_\e$ this is a basis of monomials
\[
E^{\underline  m}=E_{\beta(N)}^{m_N}\dots E_{\beta(1)}^{m_1}
\]
where $m_i\geq 0$. For $\U^-_\e$ this is a basis of monomials
\[
F^{\underline m}=F_{\beta(1)}^{m_1}\dots F_{\beta(N)}^{m_N}
\]
where $m_i\geq 0$. For $\U^0_\e$ we choose a natural
basis of Laurent monomials
for ${\underline p}=(p_1,\dots, p_r) \in {\mathbb Z}^r$ the monomial is
$K^{\underline p}=K_1^{p_1}\dots K_r^{p_r}$.

There is a linear isomorphism  $\U_\e\simeq \U^-_\e\otimes
\U^0_\e \otimes \U^+_\e$, the PBW basis in $\U_\e$ is the tensor product of
bases described above.\smallskip

   Let moreover consider $\B^+_e:=\U^0_\e
\otimes \U^+_\e,\ \B^-_e:=\U^0_\e \otimes \U^-_\e.$\smallskip

By the defining relations it follows that $\U^+_e,\U^-_\e,\B^+_e ,\
\B^-_e$  are subalgebras while $\B^+_e ,\ \B^-_e$ are even
sub-Hopf algebras. It is known that the subalgebras
\begin{itemize}
\item $Z^+_0\subset \B^+_\e$ generated by $E^\ell_\alpha, K_i^\ell$
\item $Z^-_0\subset \B^-_\e$ generated by $F^\ell_\alpha, K_i^\ell$
\item $Z_0\subset \U_\e$ generated by $E^\ell_\alpha$,
$F^\ell_\alpha$ and $K^\ell_i$,
\end{itemize}
are central and are Hopf subalgebras \cite{D0}, we will recall to
which groups these Hopf algebras correspond.

The algebras $\B^+_\e$, $\B^-_\e$ and $\U_\e$ are
CH-Hopf algebras with trace subalgebras $Z^+_0$, $Z^-_0$ and $Z_0$
respectively,
and are free of respective ranks  $\ell^{N+r}, \ell^{N+r}, \ell^{2N+r}$
over their trace subalgebras \cite{D0}.

\subsection{Structure of the center}

Let $G$ be the simply connected group associated to $\mathfrak g$, $T$
a maximal torus of $G$ and $W$ its Weyl group. $U^+,U^-$
the unipotent radicals of opposite Borel subgroups $B^+, B^-$.

Let us recall \cite{D1} that, as an Hopf algebra, $Z_0$ is the
coordinate ring of the
{\it dual group} $H$ which is the subgroup of $B^+\times B^-$  kernel
of the composed
homomorphisms $B^+\times B^-\stackrel{\mu}\longrightarrow T\times
T\stackrel{m}\longrightarrow T$ where $m$ is multiplication and
$\mu$ the quotient modulo the unipotent radical. As a variety, $H$ is
identified to $U^-\times T\times U^+$

Furthermore  $Z_0^+,Z_0^-$ are the coordinate rings of the two quotients
$B^-=H/U^+,$ $ B^+=H/U^-$.

The center $Z_\e$ of  $\U_\e$ is described in  \cite{D1}, the one of
$\U^+_\e$ will be presented later in
this paper.

Let us recall briefly the description of  $Z_\e$.

$Z_\e$ contains also another subalgebra  $Z_1$ (specialization of the central
elements for the generic value of $q$).  $Z_1$ is  identified to the
coordinate ring
of the  quotient $T/W$ (isomorphic to $G//G$ the quotient under
adjoint action).

The morphism  $B^+\times B^-\to G,\ (x,y)\to xy^{-1},$   restricted
to $H$ induces an
\'etale covering  $\rho:H\to G$, in coordinates $\rho: U^-\times
T\times U^+\to G$
is given  by
$\rho(u,t,v):=ut^2v^{-1}$.

Let $Z':=Z_0\cap Z_1$.  In \cite{D1} it
is proved that $Z_\e=Z_0\otimes_{Z'}Z_1$, moreover there is
the following geometric interpretation of this tensor product.

    a) $Z'$ is identified to the coordinate ring of the
quotient $T/W$ (isomorphic to $G//G$), under the composite map:
$$H\stackrel{\rho}\longrightarrow G\stackrel{\pi}\longrightarrow G//G=T/W $$
b)
The
   $\ell$ power map
$t\to t^\ell$ factors to the quotient giving a map $\ell:T/W\to T/W$
which at the level of coordinate rings induces the
inclusion $Z'\subset Z_1$.

c) From this we get that   $X:=V(Z_\e)$ is
    the schematic
fiber product:
$$\begin{array}{ccr}X&^{\stackrel{p}{\rule{40pt}{0.4pt}}
\!_{\longrightarrow}} & T/W\\
\sigma\downarrow & &\ell\downarrow\quad \\
H&\stackrel{\rho}\longrightarrow &G\stackrel{\pi}\longrightarrow T/W.
\end{array}$$

\vskip20pt

\subsection{Center of $\B^+_e$}

\bigskip

In \cite{D4}  and \cite{D3} (where in fact more general algebras
are studied), it is proved that the degree of the algebra
$\U^+_\e$ (resp.
$\B^+_\e$ )  over the respective center, is
$\ell^{|\Delta^+|-s\over 2}$ (resp.
$\ell^{|\Delta^+|+s\over 2}$), where
$s$ is the number of orbits of the permutation $-w_0$  on the set
$\Delta$ of simple roots, and $w_0$ is the longest element of the
Weyl group.  In particular, since the dimension of $\B^+_\e$  over
$Z^+_0$ is $\ell^{|\Delta^+|+n }$  where $n$ is the rank of the
group, we have that the center $Z^+$ of $\B^+_\e$ has dimension
$\ell^{n-s}$ over $Z^+_0$.

For types different from $A_n,\ D_n,\
E_6$ we have $w_0=-1$ and so  $Z^+=Z^+_0$.

     Otherwise there is a bigger center
which we want to describe, for type $A_{2m}$  we
have $s=m$ and for $A_{2m+1}$ we have $s=m+1$, for $D_n$ we have
$s=n-1$  and for $E_6$ we have $s=4$. \bigskip

In order to compute the center of $\B^+_e$ we need to identify this
algebra with the
so called {\it quantized function algebra} $F_\e[B^{ -}]$.
\smallskip

The construction of a function algebra is a general construction on
Hopf algebras.

Given a Hopf algebra  $\mathcal H$ and a class of finite dimensional
representations
closed under direct sum and tensor products one considers the space
$\hat {\mathcal H}$, of
linear  functions on
$\mathcal H$ spanned by the {\it matrix coefficients} $c_{\phi,v}$.

Here $v$ is a
vector in a representation $V$ and $\phi\in V^*$, the function $c_{\phi,v}$ is
defined by:
$$c_{\phi,v}(h):=\langle\phi\,|\,hv\rangle. $$
$\hat { \mathcal H}$ is also a Hopf algebra, dual to $ \mathcal H$ and
called {\it function algebra}.\medskip

In \cite{DL}  this theory is developed for the algebras of $\B^a_q,\
\B^{a,-}_q$ first at
$q$ generic. One obtains the algebras $F_q[G]$ and $F_q[B^-]$ which can be
specialized to $F_{\e}[G]$ and $F_{\e}[B^-]$ when $q$ is specialized
to a primitive
$\ell$-th root of unity $\e$.

Remark that, as $\B_\e$ is an Hopf subalgebra of $\U_\e$ so
$F_{\e}[B^-]$ is a quotient Hopf algebra of  $F_{\e}[G]$ (the same
holds for $q$
generic).

   In
particular for every dominant weight
$\lambda$ one has in $F_q[G]$ and in $F_{\e}[G]$, the matrix
coefficients of the Weyl
modules
$V_\lambda$ for the Lusztig divided power form of
$\U^a_q$.
\bigskip

An important ingredient is Drinfeld's duality which gives the
following canonical pairing
between $\B^{a,-}_q$ and $\B^{s,+}_q$ :
$$(\prod_{i=N}^1 F_{\beta_i}^{h_i}K_\alpha,
\prod_{i=N}^1 E_{\beta_i}^{h_i}K_\beta)=q^{ -(\alpha|\beta)}\prod
(h_i)_{q_{\beta_i}^{2}}!(q_{\beta_i}^{-1}-q_{\beta_i})^{-h_i} $$
(where $q_{\beta}=q^{{(\beta,\beta)}\over 2}$ and
$(h)_q={{q^h-1}\over {q-1}}$) and 0 otherwise.

  From \cite{DL} we have:
\begin{theorem}
The algebra $\B_q^{s,+} $ under this pairing is identified
   to  $F_q[B^-]$.

These isomorphisms specialize at $q$ a root of 1 giving an isomorphism between
$\B_\e^+$ and $F_\e[B^-]$.

\end{theorem}

In view of this theorem we compute the center of $F_\e[B^-]$.

We start from some identities at $q$ generic.

   From the theory of the $R-$matrix one has an immediate implication on the
commutation rules among the elements
$c_{\phi,v}$. Assume that $v,w$ have weights $\mu_1,\mu_2$ and that
$\phi,\psi$ have weights $\nu_1,\nu_2$ with respect to the
action of the elements $K_i$ \cite{LS}. Then:
$$c_{\phi,v}c_{\psi,w}=
q^{-(\mu_1|\mu_2)+(\nu_1|\nu_2)}c_{\psi,w}c_{\phi,v}+\sum
c_{\psi_i,w_i}c_{\phi_i,v_i}$$ where
$$\psi_i\otimes \phi_i=p_i(q) ( M_i(E)\otimes M_i(F))\psi\otimes
\phi,\  w_i\otimes v_i=p'_i(q) ( M'_i(E)\otimes M'_i(F))w\otimes
v $$ where the $p_i,p'_i$ are in $\mathbb  C(q)$  and $M_i,M_i'$ are
monomials of which at least one is not constant.

For each dominant weight $\lambda$ we have an irreducible
representation $V_\lambda$, we   choose  for each $\lambda$   a
highest
weight vector $v_\lambda$. We make the convention that $\phi_\lambda$
denotes a dual vector, so it is a lowest weight vector in
the dual space and it has weight $-\lambda$.

   Take
$\phi_{w_0\lambda}$ dual of a vector
$v_{w_0\lambda}$. In the commutation take
$c_{\phi_{w_0\lambda},v_\lambda}$ and a matrix coefficient
$c_{\phi,v}$ where $\phi$
has weight
$\nu$ and $v$ has weight $\mu$:

$$c_{\phi,v}c_{\phi_{w_0\lambda},v_\lambda}= q^{-(\mu |\lambda)-(\nu
|w_0\lambda)}c_{\phi_{w_0\lambda},v_\lambda}c_{\phi,v} $$
Set $\Delta_\lambda:= d_{\phi_{w_0\lambda},v_\lambda}$ and notice
that, from the previous formula we have that;
$$\Delta_\mu\Delta_\lambda=q^{-(\mu |\lambda)-(-w_0\mu
|w_0\lambda)}\Delta_\lambda\Delta_\mu=\Delta_\lambda\Delta_\mu $$
\begin{lemma}
   $\Delta_\lambda\Delta_\mu=k\Delta_{\lambda+\mu}$, $k$ a constant.
\end{lemma}
   \proof  By definition of multiplication between matrix coefficients:
$$c_{\phi_{w_0\lambda},v_\lambda}c_{\phi_{w_0\mu},v_\mu}=c_{\phi_{w_0\lambda}\otimes 
\phi_{w_0\mu},v_\lambda\otimes v_\mu}$$ Now
in the representation $V_\lambda\otimes V_\mu$ the highest weight
vector $v_\lambda\otimes v_\mu$ generates the irreducible module
$V_{\lambda+\mu}$  and in the dual we have a similar picture, thus the
matrix coefficient is only relative to this submodule.\qed

We can thus normalize the choices of the $v_\lambda$ so that
$$\Delta_\lambda\Delta_\mu= \Delta_{\lambda+\mu}$$

   \begin{lemma} A matrix coefficient $c_{\phi,v}$ where $\phi$ has weight
$\nu$ and $v$ has weight $\mu$ vanishes on $\B^{a,-}_q$ if
$\nu\not\geq -\mu$ (in the dominant order).
\end{lemma}
   \proof A form of weight $\nu$ vanishes on all vectors which do not
have weight $-\nu$, the vectors of $\B^{a,-}_qv$  have weights
$\leq \mu$ in the dominant order, so the matrix coefficient is 0
unless $-\nu\leq \mu$  or $\nu\geq -\mu$.\qed

Let us now denote by $d_{\phi,v}$  the restriction of
$c_{\phi,v}$  as function on  $\B^{a,-}_q$.

Take now $d_{\phi_\lambda,v_\lambda}$  as function
$$d_{\phi_\lambda,v_\lambda}( \prod_{i=N}^1
F_{\beta_i}^{h_i}K_\alpha)=<\phi_\lambda,\prod_{i=N}^1
F_{\beta_i}^{h_i}K_\alpha
v_\lambda>=\begin{cases} 0\quad\text{if}\quad \sum h_i>0\\ q^{<\alpha,\lambda>}
\end{cases}$$

   So, under the canonical pairing we have the identification
$$d_{\phi_\lambda,v_\lambda}=K_{-\lambda},$$

   \begin{lemma} A matrix coefficient $d_{\phi,v}$ where $\phi$ has weight
$\nu$ and $v$ has weight $\mu$ and
$\nu \geq -\mu$ (in the dominant order) is identified to a linear
combination of elements $\prod_{i=N}^1
E_{\beta_i}^{h_i}K_\alpha$ where $\sum_i h_i\beta_i=\mu+\nu$.
\end{lemma}

   \proof We have that $<\phi ,\prod_{i=N}^1 F_{\beta_i}^{h_i}K_\alpha
v >$ is 0 unless $ \mu-\sum_ih_i\beta_i=-\nu$ therefore in
the duality 1.3.3 only the terms described can occur.\qed

We have $$K_\lambda \prod_{i=N}^1 E_{\beta_i}^{h_i}K_\alpha=q^{
(\lambda,\sum_ih_i\beta_i)}\prod_{i=N}^1 E_{\beta_i}^{h_i}K_\alpha K_\lambda $$
\smallskip

Therefore, from the previous Lemma,  we get:
$$d_{\phi_\lambda,v_\lambda}d_{\phi,v}= q^{
-(\lambda,\mu+\nu)}d_{\phi,v}d_{\phi_\lambda,v_\lambda},\quad
K_\lambda d_{\phi,v}=
q^{
   (\lambda,\mu+\nu)}d_{\phi,v}K_\lambda $$ Set   $T_\lambda:=\Delta_\lambda
K_{-\lambda}$, from the previous commutation relations we get:

$$d_{\phi,v}T_\lambda = q^{-(\mu |\lambda)-(\nu
|w_0\lambda)+(\lambda,\mu+\nu)}T_\lambda d_{\phi,v}=q^{  (\lambda-w_0\lambda,
\nu)}T_\lambda d_{\phi,v} $$

   From the previous relations and the fact that $\Delta_\lambda
\Delta_\mu=\Delta_{\lambda+\mu}$  and  commute we have
   \begin{proposition}
$$ T_\lambda T_\mu= q^{(\lambda, w_0\mu- \mu)}T_{\lambda+\mu},\ $$
\end{proposition}

We can now introduce the elements

$$A_{h,\lambda}:=T_\lambda^hT_{-w_0\lambda}^{l-h}=T_{h \lambda
   +( l-h)(- w_0\lambda)}$$
    and compute the commutation relations with a matrix coefficient
$d_{\phi,v}$ where $\phi$ has weight
$\nu$ and $v$ has weight $\mu$:

$$d_{\phi,v}A_{h,\lambda}=q^{  h(\lambda-w_0\lambda, \nu)}q^{  (\ell
-h)(-w_0\lambda+\lambda, \nu)} A_{h,\lambda} d_{\phi,v}=q^\ell A_{h,\lambda}
d_{\phi,v} $$
\begin{proposition}If we specialize $q$ to an $\ell$ root of 1 we obtain that
    $A_{h,\lambda} $ is in the center.
\end{proposition}
Remark that, if $\lambda=-w_0\lambda$ we have  $ T_{h \lambda
   +( l-h)(- w_0\lambda)}=  T_{l \lambda
    }\in Z_0$ and notice that, since $\ell$ is odd
$T_{\ell\lambda}=T_{\lambda}^\ell$.\bigskip

To understand $\Delta_{\ell\lambda},$\quad (and also $T_{\ell
\lambda}=\Delta_{\ell \lambda} K_{-\ell\lambda}$) we must use the
Frobenius isomorphism, so that we identify this element to the
classical matrix  coefficient for
$\lambda$ which we will denote by  $\delta_{
\lambda}=C_{\phi_{w_0\lambda},v_\lambda}$  (denote by $C$ the
classical matrix
coefficients), by abuse of notations we denote by the same symbols
the vectors and forms in the classical representation.

Now recall that, for an algebraic group $G$, the function algebra
$\mathbb  C[G]$ has a left and a right $G$ action which in
terms of functions or of matrix coefficients are
$$(h,k)f(g):=f(h^{-1} gk),\quad (h,k)c_{\phi,v}=c_{h\phi,k v} $$
when $G$ is semisimple and simply connected, we can
exponentiate the action and identify the function algebra on
$U(\mathfrak  g)$ with the function algebra on $G$. For every dominant
weight
$\lambda$ we have an irreducible representation
$V_\lambda$ and the embedding given by matrix coefficients:
$$i_\lambda:V_\lambda^*\otimes V_\lambda\to \mathbb  C[G] $$
    The element   $c_\lambda:=C_{\phi_{w_0\lambda},v_\lambda}$, respect
to the left and right actions of $B^+\times B^+$ is an
eigenvector of weight
$-w_0(\lambda),\lambda$. In particular we can analyze it for the
fundamental weights. The following is well known
   \begin{proposition}

   The elements $c_{\omega_i}$ for the various fundamental weights are
irreducible elements whose divisors are the closures of the
codimension 1 Bruhat cells of $G$.
\end{proposition}
   \proof Let us recall one possible proof for completeness.  The ring
$\mathbb  C[G]$ is a unique factorization domain (cf. ),
the elements that are
$B^+\times B^+$   eigenvectors will then factor into irreducible
$B^+\times B^+$   eigenvectors. But these elements coincide up
to constant with the elements
$\delta_{ \lambda}$ hence the first statement is due to the fact that
the fundamental weights are free generators of the monoid
of dominant weights.

For the second part we have exactly $n=rk(G)$  codimension 1 Bruhat
cells of $G$ which must have
equations which are $B^+\times B^+$ eigenvectors.  In fact one can
identify more precisely the
correspondence (cf. ).

Another interpretation is with the Borel Weil theorem and identifying
the $c_\lambda$ with sections of line bundles on the flag
variety.\qed

Now when we restrict to $B^-$ we can exploit the fact that
$B^- U^+$ is open in $G$,  functions invariant under right $U^+$
action are identified to functions on $B^-$ we deduce that also
   \begin{proposition}\label {fund} The elements $d_{\omega_i}$,
restriction to $B^-$ of the elements
$c_{\omega_i}$,
$i=1,\ldots, n$, are irreducible elements whose divisors are the
closures of the codimension 1 Bruhat cells of
$G$ intersected with
$B^-$.

The restriction    $t_\lambda(u,t)$ to $B^-=U^-\times T$ of $T_{\ell
\lambda}=\Delta_{\ell \lambda} K_{-\ell\lambda}$ is a function only 
of $u$ and independent of $t$.
\end{proposition}
\proof  We have already proved the first part, for the second remark 
that, by definition
$\Delta_{\ell \lambda}$ transforms under right action of $T$ through 
the character $\chi_\lambda$ so
$\delta_\lambda(t,u)=g_\lambda(u)\chi_\lambda(t)$, but 
$K_{\ell\lambda}$ restricts to the character
$\chi_\lambda(t)$ hence the claim.\vskip10pt

\noindent {\bf Example.} For $SL(n)$  the fundamental representation
$\wedge^iV$ the highest weight
vector
$e_1\wedge e_2\wedge\dots\wedge e_i$ the lowest weight vector
$e_{n-i+1}\wedge e_{n-i+2}\wedge\dots\wedge e_n$  the matrix
coefficient is the determinant formed in the triangular matrix
$(x_{ij})$, by the determinant of the first $i$ rows and the last $i$ columns:
$$\vmatrix  x_{1,n-i+1}&\dots&x_{1,n}\\  x_{2,n-i+1}&\dots&x_{2,n}\\
\dots&\dots&\dots\\x_{i,n-i+1}&\dots&x_{i,n}\\ \endvmatrix
$$

\vskip30pt

   \begin{lemma}\label{bd}  Let $A$ be a unique factorization domain and a
Cohen Macaulay ring of characteristic 0 (or prime to $\ell$) and
containing the $\ell$ roots of 1.

Let
$f_1,f_2,\dots,f_k\in A$ be distinct irreducible elements and
$$R:= A[t_1,\dots,t_k]/(t_1^\ell-f_1,
t_2^\ell-f_2,\dots,t_k^\ell-f_k)$$ then $R$ is
a normal domain, Galois extension of
$A$ with Galois group $\mathbb  Z/(\ell)^k$.

\end{lemma}
   \proof Clearly $R$ is free over $ A$ of rank $\ell^k$ and $\mathbb
Z/(\ell)^k$
acts as symmetry group.  We need only show that $R$ is a normal domain.

1) First of all
$R$ is a complete intersection hence it is Cohen Macaulay.

2) Next we will prove that
it is smooth in codimension 1 which will prove that it is a normal ring.

3) Finally we
prove that its spectrum is connected which will imply that it is a domain.

2) Since $A$ is normal we can restrict our analysis to the smooth
locus and choose a
regular system of parameters $x_1,\dots,x_m$. Consider the
Jacobian matrix (e.g.
$k=3$):
$$\vmatrix \ell t_1^{\ell-1}&0&0& \pd{f_1}{x_1} & \pd{f_1}{x_2}
&\dots & \pd{f_1}{x_m}  \\ 0&\ell t_2^{\ell-1}&0 & \pd{f_2}{x_1}
&
\pd{f_2}{x_2} &\dots & \pd{f_2}{x_m}     \\ 0&0&\ell
t_3^{\ell-1}&\pd{f_3}{x_1} & \pd{f_3}{x_2} &\dots &
\pd{f_3}{x_m}
\endvmatrix$$ On the open set where the $f_i$ are non zero the first
determinant is not zero and thus $R$ is smooth in
codimension 0 and hence reduced. In a smooth point of the subvariety
$f_i=0$ where
also $\prod_{j\neq i}f_j\neq 0$, we also have a non zero maximal
determinant so also
these points are smooth. The complement has codimension at least 2.

3) Finally we have to prove connectedness. Let $\overline F$ be an
algebraic closure of $F$,  we can argue as follows, let
us consider the ring $\overline R\subset
\overline F$ obtained from $R$ by adding $\ell$ roots $b_i$ of $f_i$
we have clearly a homomorphism of $R$ onto $\overline R$.
Let $Q,G$ be the quotient fields of $R,\overline R$, clearly $G$ is a
Galois extension of $Q$ with Galois group a subgroup
$\Gamma$  of $\mathbb  Z/(\ell)^k$, it is clearly enough to show that
this subgroup is $\mathbb  Z/(\ell)^k$ itself. \smallskip

Let $\mathcal  M:=\{b_1^{h_1}\dots b_k^{h_k},  (h_1,h_2,\dots,h_k)\in
\mathbb  Z^k\}$ and $\epsilon=e^{2\pi i\over \ell}$, we
identify $\mathbb  Z/(\ell) $ with the multiplicative group generated
by $\epsilon$ and  have a pairing:
$$\Gamma\times \mathcal  M\stackrel{p}\longrightarrow \mathbb  Z/(\ell),\quad
p(\sigma,M):=\sigma(M)M^{-1}
$$ This pairing factors through $\mathcal  M^\ell$ and if by contradiction
$|\Gamma|<\ell^k$ we must have an element $M=b_1^{h_1}\dots b_k^{h_k},\ 0\leq
h_i<\ell$  not all the $h_i=0$ which is in the kernel of the pairing, hence by
Galois theory (and the fact that the element is integral over
$A$) $M\in A$. Then
$$M^\ell=f_1^{h_1}\dots f_k^{h_k}. $$ Factoring $M$ into
irreducibles, this implies
that $\ell\,|\,h_i$ for all $i$   a contradiction unless all $ h_i=0$.\qed

Let $a_i:=T_{\ell \omega_i}\in Z_0^+$, $i=1,\ldots ,n$. Notice that,
under the identification of
$Z_0^+$ with $\mathbb C[B^-]$, the elements $a_i$ coincide up to a
unit with the elements $d_{\omega_i}$
defined in Proposition \ref{fund},  so that their divisors are
irreducible and distinct. Consider the
algebra
$R:=Z_0[b_1,\dots,b_n]$ with $b_i^\ell=a_i$. Lemma \ref{bd} implies that
$R$ is a normal domain on which acts the Galois group $\mathbb
Z/(\ell)^n$. Let $\tau:\omega_i\to
-w_0\omega_i$ be the standard involution of fundamental weights.
$\tau$ induces  an involution of
the factors of $\mathbb  Z/(\ell)^n$. Let
$\Gamma$ be the invariant subgroup, it is made of those $r-$tuples
which have the
same entry in   the orbits of $\tau$ (made of 1 or 2 elements). The
invariants under
$\Gamma$ are spanned by the monomials $b_1^{h_1}b_2^{h_2}\dots
b_n^{h_n}$ which when
$i,j$ are an orbit of $\tau$ have exponents $h_i+h_j\equiv
0,\quad\text{mod}\ \ell$.

Thus $R^\Gamma$ is isomorphic to the algebra generated by the
elements $ T_{h \lambda
   +( \ell-h)(- w_0\lambda)}$.
\begin{theorem} The center of $\B_\e^+ $ is the algebra $Z_\e^+$ 
generated by $Z_0^+$ and by  the
elements $ T_{h \lambda
   +( \ell-h)(- w_0\lambda)}$ .
\end{theorem}
   \proof The algebra $Z_\e^+$ being isomorphic to $R^\Gamma$, is
normal. Also, since the degree of
$B_\e^+
$ equals $\ell^{|\Delta^+|+s\over 2}$, where
$s$ is the number of orbits of  $\tau$, while $B_\e^+
$   has  rank $\ell^{|\Delta^+|+n}$  over
$Z_0^+,$ we deduce that the rank of $Z_\e^+$ over
$Z_0^+,$ equals the rank of the center, so $Z_\e^+$ is the center.\qed

\vskip20pt

\section{Clebsch-Gordan decompositions for generic
representations of quantized universal enveloping algebras at roots of 1}

\subsection{Compatibility for $\U^+_\e\subset \U_\e$  }

\begin{proposition}The natural map
$\U+_\epsilon\otimes_{Z^+_0}Z_\epsilon\to \U_\e$ is injective.
\end{proposition}

\proof Recall that $\U^+_\e$ is free over $Z^+_0$. As a linear basis in
$\U^+_\e$ over $Z^+_0$ we can choose PBW elements $b=\prod_{\alpha\in \Delta_+}
E_\alpha^{m_\alpha}$ with $0\leq m_\alpha<\ell$. The center $Z_\e$ is
generated by $Z_1$ and $Z_0$ where $Z_1$ is the ``specialization at $q=\e$
of the center for generic $q$. Any element $z\in Z_1$ is completely
determined by its component $\phi_{0,0}\in \U^0_\e$ and it is of the form

\begin{equation}\label{Z1}
z=\phi_{0,0}+\sum_{{\underline r},{\underline k}}
E^{\underline r}\phi_{{\underline r},{\underline k}}F^{\underline k},
\ \ \phi_{{\underline r},{\underline k}}\in \U^0_\e
\end{equation}

Moreover $Z_\e=Z_0\otimes_{Z_0\cap Z_1}Z_1$.

We want to prove that if $\sum bz_b=0$ where $b$ are elements of the
PBW basis and $z_b\in Z_\e$, then for any $b$, $z_b=0$.

%More or less the idea is as follows:
%\begin{enumerate}
%\item we know $\U^+_0$ is free over $Z^+_0$,
%\item the proposition would be easy to prove if the coefficients $b$
%were in $Z_0$ instead
%of $Z_\e$,
%\item  apply the coproduct to $\sum bz_b$, and isolate the part in
%$\U^+_\e \otimes\U^0 \U^-_\e$. From the
%general form (\ref{Z1}) of elements of
%$Z_1$ , we observe that the ``higher degree part'' $\sum E^{\underline r}
%\phi_{{\underline r},{\underline k}}F^{\underline k}$ will add
%only ``higher terms'' to the left side of the $\otimes$ after applying
%the comultiplication.
%\end{enumerate}

%Now the proof:
%\begin{itemize}

    Recall that we have a convex ordering on $\Delta_+$. In the product
defining PBW elements we choose the decreasing order of
$E_\alpha^{m_\alpha}$. This provides a total ordering on PBW elements defined
by the lexicographic ordering of $E_\alpha$.

%\item Apply the coproduct $\Delta$ to $\sum bz_b$ and pick the component
%in $\U^+_\e\otimes \U^0_\e\U^-_\e$.

\begin{lemma}\label{PBW}
Let $b=\prod_\alpha E_\alpha^{m_\alpha}$ be a PBW element and
$\beta\in \Delta_+$ ; then $bE_\beta$ is a linear combination of
PBW elements which are  greater than $b$.
\end{lemma}
\proof We have $b=E_{\beta(1)}^{m_1}\dots E_{\beta(N)}^{m_N}$ where
$\beta(1)>\dots >\beta(N)$ in the convex ordering of $\Delta_+$. Then
\begin{itemize}
\item if $\beta(N)\geq \beta$ the statement is clear,
\item if $\beta(N)<\beta$ we use the commutation relation

\begin{eqnarray}
E_{\beta(N)}E_\beta & = & q^{(\beta(N),\beta)}E_\beta E_{\beta(N)} \\
+ &\sum_{\beta(N)<
\gamma_1<\dots<\gamma_s<\beta}& a_{{\gamma_1,\dots, \gamma_s}{r_1,\dots,r_s}}
E_{\gamma_s}^{r_s}\dots E_{\gamma_1}^{r_1}
\end{eqnarray}
It is clear that here $E_\beta E_{\beta(N)}>E_{\beta(N)}$ and
$E_{\gamma_s}^{r_s}\dots E_{\gamma_1}^{r_1}>E_{\beta(N)}$.
Now we can iterate this process to reorder monomials and the
lemma follows.
\end{itemize}
\qed

\begin{itemize}
\item Coming back to the proof of the proposition we look at the PBW
elements $b$ for which $z_b\neq 0$ and let $b_0$ be the minimal
among them. Then let us apply the coproduct  and take the component
which belongs to $\U^+_\e\otimes \U^0_\e \U^-_\e\subset \U_\e\otimes \U_\e$.

\item Due to triangular decomposition in $\U_\e$,
each $\Delta(b)$ will contribute only by $b\otimes 1$ to this component.
The element (\ref{Z1}) will contribute as
\[
1\otimes \phi_{0,0}+\sum E^{\underline r}\otimes
\phi_{{\underline r},{\underline k}}F^{\underline k}
\]
and any element in $Z_0$ which is always a polynomial in
$x_\alpha=E_\alpha^\ell, \ y_\alpha=F_\alpha^\ell, \ z^\pm_i=k^{\pm \ell}_i$
will contribute as a polynomial in $x_\alpha\otimes 1, \ 1\otimes
y_\alpha, 1 \otimes z^\pm_i$.

The minimal term in the left side of the tensor product is of the form
$b_0P(x_\alpha)$ for some polynomial $P$ because, from lemma \ref{PBW},
the terms coming from $Z_1$ will contribute by $1$ up to bigger terms.
The contribution of
PBW elements $b>b_0$ to the left component of the tensor product
will have their minimal monomial of exactly the same form.
\end{itemize}

Now the proposition follows  from the freeness of
$\U^+_\e$ over $Z^+_0$.
\qed

\subsection{Compatibility of comultiplication}
In this paragraph we will strengthen Theorem \ref{com}
as follows: 
\begin{theorem}\label{bel}  Comultiplication $\Delta:\U_\e\to \U_\e\otimes
\U_\e$  is  compatible (with trace) when we think of
$\U_\e$ as $Z_0$ algebra and $\U_\e\otimes
\U_\e$ as
$Z_\e\otimes  Z_\e$ algebra.
\end{theorem}
Using  \ref{comp} ii) we need to show:

\begin{proposition} \label{tre}
$\Delta(Z_\e)\otimes_{\Delta(Z_0)} (Z_\e\otimes
Z_\e)$ is a normal domain.\end{proposition}

\proof
We have recalled the analysis of \cite{D4} in 6.3 and in particular
the fiber product diagram:
$$\begin{array}{ccr}X& ^{\stackrel{{\huge p}}{\rule{40pt}{0.4pt}}
\!_{\longrightarrow}} & T/W\\
\sigma\downarrow & &\ell\downarrow\quad \\
H&\stackrel{\rho}\longrightarrow &G\stackrel{\pi}\longrightarrow T/W.
\end{array}$$
where $X$ is the spectrum of $Z_\e$.

The analysis shows that one can define a regular locus in all these varieties.

$G^{reg}$ is the usual set of regular elements (i.e. elements with
conjugacy class of maximal dimension).

Finally we set $$H^{reg}:= \rho^{-1} G^{reg},\quad
X^{reg}:=\sigma^{-1} H^{reg}.$$
In the restricted fiber product diagram
$$\begin{array}{ccr}X^{reg}&^{\stackrel{p}{\rule{40pt}{0.4pt}}
\!_{\longrightarrow}} & T/W\\ 
\sigma\downarrow & &\ell\downarrow\quad \\
H^{reg}&\stackrel{\rho}\longrightarrow
&G^{reg}\stackrel{\pi}\longrightarrow T/W.  
\end{array}$$
 by \cite{St}, we know
that the subset
$G^{reg}$ of $G$ of regular elements has a complement of codimension
2 and that the map $\pi$ restricted to $G^{reg}$ is a smooth map.  It follows that
$\pi\circ \rho$ and also $p$ are smooth and since $T/W$ is smooth,  that all
varieties in this regular diagram are smooth. 
Since the complement of
$H^{reg}$ in $H$ has codimension
$\geq 2$ and $\sigma$ is finite,    the complement of $X^{reg}$ in $X$  has
also codimension
$\geq 2$.
\smallskip

   We know that
the ring $Z_\e$ is presented as a complete intersection over
$Z_0$  by the fiber product diagram and it is free of finite rank.
Thus we have that, $Z_\e$ and  $Z_\e\otimes
Z_\e$ are normal Cohen Macaulay domains. For the same reasons
$\Delta(Z_\e)\otimes_{\Delta(Z_0)} (Z_\e\otimes
Z_\e)$ is a complete intersection, hence a Cohen-Macaulay ring.

$\Delta(Z_\e)\otimes_{\Delta(Z_0)} (Z_\e\otimes
Z_\e)$ is the coordinate ring of the schematic fiber product $Y$:
$$\begin{array}{ccr}Y&^{\stackrel{p}{\rule{30pt}{0.4pt}}
\!_{\longrightarrow}} &\!\!\! \!\!\! \!\!\! \!\!\! \!\!\!
\!\!\! \!\!\! \!\!\! \!\!\! \!\!\! \!\!\! X\\
\downarrow &  &\quad\quad\sigma\downarrow \\
X\times X&\stackrel{\sigma\times \sigma}\longrightarrow &H\times
H\stackrel{m}\to H.
\end{array}$$ where $m$ is multiplication. This can also be presented
as the unique fiber product map:
$$\begin{array}{lclcl}Y&^{\stackrel{p}{\rule{30pt}{0.4pt}}
\!_{\longrightarrow}} & \qquad \quad X
& ^{\stackrel{p}{\rule{40pt}{0.4pt}} \!_{\longrightarrow}}  & T/W\\
q\downarrow &  &\quad\quad\sigma\downarrow & &\ell\downarrow\quad \\
X\times X&\stackrel{\sigma\times \sigma}\longrightarrow &H\times
H\stackrel{m}\to
H&\stackrel{\rho}\longrightarrow &G\stackrel{\pi}\longrightarrow T/W.
\end{array}$$

   We
will then apply Serre's criterion \cite{Se}, and prove that $Y$
is smooth in codimension  1, which will show that its coordinate ring
is a {\it normal
ring}.

   By the homeomorphism $\nu:H\times H\to
H\times H,\ \nu(x,y):=(x,xy)$ it follows that the open set
$\mathcal A:=\{(x,y)\in H\times H\}$ with $x,y,xy$ regular  has a
complement of codimension 2; on this set the map  $(x,y)\to \pi\circ \rho(xy)$
of $H\times H\to T/W$ is a smooth map.

On the open set $\mathcal B:=(\sigma\times\sigma)^{-1}\mathcal A$ the
composite map $\pi\circ  \rho\circ
(\sigma\times\sigma)$ is smooth, and since the map
$T/W^{reg}\stackrel{\ell}\to T/W^{reg}$ is also smooth we deduce
that $\mathcal C:=q^{-1}\mathcal B$ is smooth.  Since $q$ is finite
it follows that the complement of $\mathcal C$ in
$Y$ has codimension $\geq 2$ hence $Y$ is a normal variety
and it follows that
$\Delta(Z_\e)\otimes_{\Delta(Z_0)} (Z_\e\otimes
Z_\e)$ is a normal ring.

It remains (as in Lemma \ref{bd} ) to finish the argument and prove
that $Y$ is connected which implies that
$\Delta(Z_\e)\otimes_{\Delta(Z_0)} (Z_\e\otimes Z_\e)$ is a normal
domain.\bigskip

Since the morphism $q$ is finite it is also proper so it suffices to
find some closed subvariety
$M$ of
$X\times X$ with the property that $q^{-1}(M)$ is connected.

First of all, we claim that we have a natural embedding:

$U^+\stackrel{i^+}\to X$ so that the diagram:
\qquad\xymatrix{& X\ar[d]_{\rho} \\
     U^+\ar[ru]^{i+} \ar[r]_v & H }
\qquad   is commutative.

In fact we see that the composed map
       $U^+\to H\to G\to G//G$  is constant with value
the class $\overline 1$ of 1,
   so it can be lifted by choosing a
   point in
the fiber of $\ell^{-1}(\overline 1)$.

We now embed  $U^+\times U^-\stackrel{j}\to X\times X$ and see that
the composed map
$$U^+\times U^-\stackrel{j}\to X\times X\stackrel{\rho\times \rho}\to
H\times H\stackrel{m}\to  H\stackrel{\sigma}\to
G$$ induces the natural inclusion by multiplication $U^+ U^-\subset G$.

In the next Lemma we show that we  have a section $M\subset U^+ U^-$
for which the composed map $M\stackrel{\pi}\to
G//G$ is an isomorphism, hence $q^{-1}(M)$, i.e. $Y$ restricted to
$M$ is
$T/W$ and connected.

\begin{lemma} There is a section $M\subset U^+ U^-$ for which the
composed map $M\stackrel{\pi}\to
G//G$ is an isomorphism.

\end{lemma}

We will need in our analysis a variation of a result of Steinberg.
Let us recall his Theorem. Let $G$ be a semisimple
simply connected group.

For $\beta$ a positive root let us denote by $X_\beta=\exp(\mathbb C
e_\beta)$ the root subgroup
associated to $\beta$. If
$\alpha_1,\dots,\alpha_n$ is the set of simple positive roots denote
by $\sigma_i$ a representative in
the normalizer of the torus
$T$ of the simple reflection   $s_i$ associated to the root
$\alpha_i$. Finally let $\pi:G\to G//G=T/W$ be the quotient under
adjoint action. Define:
$$  N=X_{\alpha_1}\sigma_1X_{\alpha_2}\sigma_2\dots
X_{\alpha_n}\sigma_n$$ the theorem of Steinberg is
that $N$ is a slice of the map $\pi$ in other words under $\pi$ $N$
is isomorphic to $T/W$.

For our purposes we have to slightly change this type of slice, we
start remarking that
$$  N=X_{\beta_1} X_{\beta_2} \dots X_{\beta_n} \sigma_1
\sigma_2\dots  \sigma_n,\quad
\beta_i:=s_1s_2\dots s_{i-1}(\alpha_i)$$ next we want to show that,
provided we possibly change the
representative $\sigma_n$, we can express $$ \sigma_1 \sigma_2\dots
\sigma_n=a_+b_-c_+,\quad
a_+,c_+\in U^+,\ b_-\in U^-.
$$ For this consider the flag variety $\mathcal B$ and in it the
point $p^+$ with stabilizer $B^+$, consider
$q:=\sigma_1
\sigma_2\dots  \sigma_n p_+$ and then
$U^+q\cap U^-p_+\neq \emptyset$ (cf. ). Thus we can find  $a_+\in
U^+,\ b_-\in U^-$ with $a_+^{-1}q=b_-p_+$ hence
$b_-^{-1}a_+^{-1}\sigma_1 \sigma_2\dots
\sigma_np_+ =p_+$ hence:
$$\sigma_1 \sigma_2\dots  \sigma_n=a_+b_- c_+t,\ \quad c_+\in U^+,\
t\in T.$$ We change then $\sigma_n$
with
$\sigma_nt^{-1}$ and get
$$\sigma_1 \sigma_2\dots  \sigma_n=a_+b_- c_+ $$ Now we obtain the
new slice $M:=c_+Nc_+^{-1}$.

   The interest for us is that $$M=c_+X_{\beta_1} X_{\beta_2} \dots
X_{\beta_n} a_+b_-\subset U^+U^-.$$
as requested.\qed
\vskip10pt

There is actually a rather interesting application of the slice that we found.

We can consider $U^+U^- $ also as subset of $H$, the canonical covering
$\sigma:H\to G$ restricted to  $U^+U^- $ is a homeomorphism to the image.
Therefore we can also consider $M\subset H$, in \cite{D1} it is shown that the
preimage of  a regular orbit of $G$ is a unique symplectic leaf in $H$
while it is the union of $\ell^n$ leaves in $X$. We deduce that we have the
regular elements in $X$ and $H$ which are unions of maximal Poisson leaves and
that:
\begin{theorem}
The set $M\subset H$  is a cross section of the set of regular Poisson leaves in
$H$.

The set $\rho^{-1}(M)\subset X$ is homeomorphic to $T/W$ and  is a cross section of the set of regular Poisson leaves in
$X$. 
\end{theorem}
\vskip20pt
\subsection{Clebsch Gordan formula}
\bigskip
    We know that $\U_\e$ has no zero divisors,  as a $Z_0$-algebra
is a free $Z_0$-module of rank $\ell^{2|\Delta_+|+n}$.  By \cite{D4}
    its center $Z_\e$ is a free $Z_0$-module of rank $\ell^{ n}$. If
$Q(Z_0)$ denotes the quotient field of $Z_0$ we have that
$Q(\U_\e)=\U_\e\otimes_{Z_0}Q(Z_0)$ is a division algebra of
dimension $\ell^{2|\Delta_+| }$ over its center
$Q(Z_\e)=Z_\e\otimes_{Z_0}Q(Z_0)$.

    Therefore:
$$[\U_\e:Z_0]= \ell^{ |\Delta_+|+n},\quad [\U_\e:Z_\e]= \ell^{
|\Delta_+| },\quad [Z_\e:Z_0]= \ell^{  n} .$$

Let $V$ and $W$ be two generic irreducible representations  of $\U_\e$
of maximal dimension $m=\ell^{|\Delta_+|}$. We want to decompose the
representation $V\otimes W$ of $\U_e\otimes \U_\e$ into irreducible
representations of the subalgebra $\Delta(\U_\e)$.

We apply the methods of Theorem \ref{CG},  recalling that $Z_0$ is a
Hopf subalgebra of $\U_\e$, but   $Z_\e$ is
only a subalgebra. So, if $V=M_P,\ W=M_Q$ where $P,Q\in V(Z_\e)$  and
$\pi(P)=x\in V(Z_0),\pi(Q)=y\in  V(Z_0)$  we know by
\ref{CG} that, for generic $x,y$:
$$M_P\otimes M_Q\simeq \oplus_{R\in \pi^{-1}(xy)} M_R^{\oplus h_R^{P,Q}}$$
we want to prove in our case:

\begin{theorem}\label {brut}
The multiplicities $h_R^{P,Q},\ R\in \pi^{-1}(xy)$, are all equal to
$\ell^{|\Delta_+|-n}$.
\end{theorem}
\proof In view of Theorem \ref{rest}, in order to prove this Theorem, since by the
generic assumption $\dim M_P=\dim M_Q=\dim M_R=\ell^{|\Delta_+|}$, 
and the degree of $\pi$ is $\ell^n$, it is enough to use the compatibility proven
in  Theorem
\ref{bel}. 

\vskip20pt
\subsection{Compatibility for $\B^+_\e\subset \U_\e$}

Using the results of 6.4 we prove now:
\begin{theorem}  i) $Z_\e^+\otimes_{Z_0}Z_\e$   is a
normal domain.

ii) The inclusion $\B^+_\e\subset \U_\e$ gives compatible algebras,
where $\B^+_\e$ is thought as $Z_0$ and $\U_\e$ as
$Z_\e$ algebras.
\end{theorem}
\proof

  From the analysis leading to 6.10 we know that:
$$Z_\e^+=Z_0[b_1,\dots,b_n]^\Gamma,\quad b_i^\ell=a_{i}.$$
thus also:
$$Z_\e^+\otimes_{Z_0}Z_\e=Z_\e[b_1,\dots,b_n]^\Gamma,\quad b_i^\ell=a_{i}.$$
therefore it is enough to prove that $Z_\e[b_1,\dots,b_n]$ is a normal domain.

Let us  argue geometrically.  Let $V$ be the variety of
$Z_\e[b_1,\dots,b_n]$, a normal variety.
   $Z_\e[b_1,\dots,b_n]$ is the coordinate ring of the schematic fiber
product $S$:
$$\begin{array}{ccr}S&^{\stackrel{s}{\rule{30pt}{0.4pt}}
\!_{\longrightarrow}} &\!\!\! \!\!\! \!\!\! \!\!\! \!\!\!
\!\!\! \!\!\! \!\!\! \!\!\! \!\!\! \!\!\! V\\
t\downarrow &  &\quad\quad v\downarrow \\
X &\stackrel{\sigma }\longrightarrow &H \stackrel{p}\to B
\end{array}$$
We have that the map $t$ is finite and flat so $S$ is Cohen
Macaulay since, as we have seen, $X$ is Cohen
Macaulay.  The map
$p$ is smooth, so the composite map
$p\circ
\sigma$ is smooth in the corresponding {\it regular elements} i.e.
outside subvarieties of codimension
2.

  From  Lemma \ref{bd} we have that also $v$ is smooth outside
codimension 2 and, since $t$ is finite we deduce that
$S$ is smooth in codimension 1 hence $S$ is a normal variety.

To prove the irreducibility  of $S$ we can, reasoning as in
Proposition \ref{tre}, restrict to the
section
$U^-\subset X$, analyze the schematic fiber product diagram:
$$\begin{array}{ccr}S'&^{\stackrel{s}{\rule{30pt}{0.3pt}}\!_{\longrightarrow}}
& V\\
t\downarrow &  & v\downarrow \\
U^- &\stackrel{\sigma }\longrightarrow &H \stackrel{p}\to \ B
\end{array}$$
and show that $S'$ is irreducible. The coordinate ring $A$ of $S'$ is
obtained from the coordinate ring
$\mathbb  C[U^-]$ of
$U^-$ by adding the
   restrictions  $\overline b_i$, of the elements $b_i$ which are of
course  $\ell$-th
roots of the  restrictions of the $a{_i}$.  Now
identify
$B^-=U^-\times T$. By the definitions of the functions
$a{_i}$ we have that these functions are invariant under right action of $T$.
It follows that
$a{_i}(u,t)=g_i(u)$, hence the $g_i(u)$ are irreducible polynomials on $U^-$
defining distinct divisors and the coordinate ring $A$ of $S'$ is:
$$A=\mathbb C[U^-][\overline b_1,\dots,\overline b_n],\quad \overline
b_i^\ell=g_i(u).$$
We can then again apply Lemma \ref{bd} and deduce that $S'$ is
irreducible, concluding
the proof of the theorem.\qed

We obtain as corollary, using Theorems \ref{ress} and
\ref{rest}, the branching rules from $\U_\e$ to $\B^+_\e$. Recall
that the degree of $\U_\e$ (as
algebra over its center) is $\ell^{|\Delta_+|}$, while the degree
of $\B^+_\e$ (as algebra over
$Z_0^+$) is $\ell^{{|\Delta_+|-s\over 2}+n}$, so the factor $r$ of
\ref{ress} is
$\ell^{{|\Delta_+|+s\over 2}}$
\begin{theorem}
Let $M$ be a semisimple trace representation of  $\U_\e$ with central
character  $\chi$ on $ Z_e$. Let $p$
the point in
$B^-$ induced by $\chi$.

i)  The restriction  of $M$ to  $\B_\e^+$ has as
associated semisimple
representation the one $N_p$  of character $p$ on the coordinate ring
$Z_0^+$ of $B^-$ with multiplicity
$\ell^{{|\Delta_+|+s\over 2}-n}$.

 ii)   For generic $p$, $N_p$ is the direct sum of the irreducible trace representation
for the center of $\B_\e^+$ with
central characters all the $\ell^{n-s}$ central characters which restricted to
$Z_0^+$ give $p$.
\end{theorem}

In particular consider a generic point $\chi\in X$. We have  a unique
irreducible
representation of
$\U_\e$ with central character  $\chi$.  When we restrict it to
$\B_\e^+$ we have a direct sum
of all the   $\ell^{n-s }$ irreducible representations, each of dimension 
$\ell^{{|\Delta_+|+s\over 2}}$, which on
$Z_0^+$ have central character $p$ and each with multiplicity
$\ell^{{|\Delta_+|+s\over 2}-n}$.\vskip10pt
\subsection{Compatibility of comultiplication for  $\B_\e^+$}  We can repeat for
$\B_\e^+$ the same analysis done in
\S 7.2 for
   $\U_\e$. Recall that:
$$[\B_\e^+:Z_0^+]=\ell^{{|\Delta_+|-s\over 2}+n},\quad 
[\B_\e^+:Z_\e^+]= \ell^{{|\Delta_+|+s\over
2}},\quad [Z_\e^+:Z_0^+]= \ell^{  n-s} $$

Let $V$ and $W$ be two generic irreducible representations  of $\B_\e^+$
of maximal dimension $m=\ell^{{|\Delta_+|+s\over
2}}$. We want to decompose the
representation $V\otimes W$ of $\B_\e^+\otimes \B_\e^+$ into irreducible
representations of the subalgebra $\Delta(\B_\e^+)$.

We apply the methods of Theorem \ref{CG},  recalling that $Z_0^+$ is a
Hopf subalgebra of $\B_\e^+$, but   $Z_\e^+$ is
only a subalgebra. So, if $V=M_P,\ W=M_Q$ where $P,Q\in V(Z_\e^+)$  and
$\pi(P)=x\in V(Z_0),\pi(Q)=y\in  V(Z_0^+)$  we know by
\ref{CG} that, for generic $x,y$:
$$M_P\otimes M_Q\simeq \oplus_{R\in \pi^{-1}(xy)} M_R^{\oplus h_R^{P,Q}}$$
we want to prove in our case:

\begin{theorem}
The multiplicities $h_R^{P,Q},\ R\in \pi^{-1}(xy)$, are all equal to
$\ell^{{|\Delta_+|+s\over
2}-n+s}$.
\end{theorem}
\proof In order to prove this Theorem, since by the generic
assumption $\dim M_P=\dim M_Q=\dim M_R=\ell^{{|\Delta_+|+s\over
2}}$, it is
enough,
by  Proposition
\ref{rest} to show the stronger statement that the inclusion
$\Delta(\B_\e^+)  \subset  \B_\e^+\otimes \B_\e^+ $  is compatible, when we
consider $\B_\e^+$ as $Z_0$ algebra but $\U_e\otimes \U_e $ as
$Z^+_\e\otimes Z^+_\e$  algebra.  Hence we need to show,   looking at the
following restriction diagram~:
\[
\begin{array}{lcl}
\Delta(\B_\e^+) &\subset &\B_\e^+\otimes \B_\e^+  \\
\cup & & \cup  \\
\Delta(Z^+_0) & \subset & Z^+_\e\otimes Z^+_\e
\end{array}
\]
     that $\Delta(\B_\e^+)\otimes_{\Delta(Z^+_0)} (Z^+_\e\otimes
Z^+_\e)$
embeds in $\B_\e^+\otimes \B_\e^+$.

Recalling Theorem 5.12 and definition 5.7, we know that
comultiplication is a compatible inclusion when we think of
$\B_\e^+$ as $Z^+_0$ algebra and $\B_\e^+\otimes \B_\e^+$ as $Z^+_0\otimes
Z^+_0$ algebras but we need the  stronger   fact that
\begin{theorem} Comultiplication $\Delta:\B_\e^+\to \B_\e^+\otimes
\B_\e^+$it is  compatible (with trace) when we think of
$\B_\e^+$ as $Z^+_0$ algebra and $\B_\e^+\otimes
\B_\e^+$ as
$Z^+_\e\otimes  Z^+_\e$ algebra.
\end{theorem}
Using  \ref{comp} ii) we need to show:

\begin{proposition} \label{trre}
$\Delta(Z^+_\e)\otimes_{\Delta(Z^+_0)} (Z^+_\e\otimes
Z^+_\e)$ is a normal domain.\end{proposition}
\proof The proof is similar to the one  of Proposition \ref{tre}. We
identify  this ring to a ring of invariants
of a ring obtained by extracting roots and prove the usual fiber 
product smoothness condition.
\vskip20pt
\subsection{The degrees and dimensions of cosets}
Here we will present hueristic arguments in favor of the idea
that, the formulae for degrees that we obtained
imply the existence of birational Darboux coordinates on the
corresponding cosets.

\subsubsection{} Let $\M_{2d}$ be a compact symplectic manifold. Geometric
quantization produces a sequence of vector spaces
$\{V_n\}_n, \ n=1,2,\dots$. The corresponding sequence $\{ End(V_n)\}$
of matrix algebras can be regarded as a quantization of the Poisson
algebra of functions on $\M_{2d}$. For large $n$ the dimension of
$V_n$ have the following asymptotic behavior
\begin{equation}\label{asym}
dim(V_n)=Vol(\M_{2d})n^d (1+O(1/n)) \ ,
\end{equation}
where $Vol(\M_{2d})$ is the symplectic volume of the symplectic
manifold.

Let $\M_{2d}=\T^{2d}$ be the $2d$-dimensional torus with coordinates
$t_1,\dots, t_{2d}\in \C, \ |t_i|=1$. Assume that the symplectic structure
on this manifold is constant:
\begin{equation}\label{sform}
\omega=\sum_{a,b=1}^{2d} \omega^{ab}\frac{dt_a}{t_a}\wedge
\frac{dt_b}{t_b}
\end{equation}
where $(\omega^{ab})$ is an integral matrix invertible over $\mathbb Z$.
Geometric quantization of this manifold produces the sequence
of vector spaces $V_n$ with $dim(V_n)=n^d$. Because the
symplectic structure is constant, the asymptotic formula
(\ref{asym}) becomes exact. One can argue that tori are typical
manifolds for which this takes place.

The complexification of $(\T^{2d},\omega)$ is the complex torus
$({\mathbb C}^*)^{2d}$ with complex holomorphic symplectic form
(\ref{sform}). The algebra of Laurent polynomials in $t_i$ is a Poisson
algebra with the brackets
\[
\{t_i,t_j\}=\omega_{ij}t_it_j
\]
where $\omega_{ij}$ is the matrix inverse to $\omega^{ij}$.

Let $q$ be a nonzero complex number. Define the algebra
$C_q(({\mathbb C}^*)^{2d})$ generated by $t_i^{\pm 1}$ with
defining relations
\[
t_it_i^{-1}=1, t_it_j=q^{\omega_{ij}}t_jt_i
\]
This family of algebras is a deformation quantization of
the Poisson algebra of functions on $({\mathbb C}^*)^{2d}$.

For a primitive root of unity $\e$ of degree $\ell$
consider the specialization of the algebra $C_\e(({\mathbb C}^*)^{2d})$
to $q=\e$. It is clear that Laurent polynomials in $t_i^{\ell}$ are in the
center of this algebra. Moreover, it is well known that they generate the
center and that $C_\e(({\mathbb C}^*)^{2d})$ is a Cayley-Hamilton algebra
over its center of degree $\ell^{d}$.

It is remarkable that the degree in this case coincides with the
square of the dimension of the space obtained by geometrical
quantization with $n=\ell$ and that this dimension coincides
with its asymptotic (\ref{asym}).

There are many examples of algebraic symplectic varieties which are
birationally equivalent to a complex symplectic torus.
It is natural to expect that if a sequence of
Cayley-Hamilton algebras quantizes such symplectic variety in a certain regular
way then degrees of such algebras will be $n^{d}$ where $d$ is
half the dimension of the complex torus. Conversely, if there is
a sequence of such Cayley-Hamilton algebras quantizing a Poisson
variety then, one can take it as an indication that the
variety is birationally equivalent to a symplectic torus.

\subsubsection{} The dimensions of multiplicity spaces which we studied all
have the form $\ell^d$ for some $d$. Hueristic arguments presented
above suggest the following conjectures
concerning the multiplicities. Let $G$ be as before, a semisimple group of rank $n$.

{\bf 1}. Let ${\mathcal O}\subset G$ be a generic conjugation
$G$-orbit in $G$ and ${\mathcal O}_-\subset B_-$ be a generic
dressing orbit of $B_+$ in
$B_-$ for the standard Poisson structure on $B_-$.
For generic orbits we have $dim({\mathcal O})=2|\Delta_+|$ and
$dim({\mathcal O}_+)=|\Delta_+|-dim(ker(w_0-id))$ where $w_0\in W$ is the
longest element of the Weyl group.

Consider the Poisson structure on $G$ which comes from $G^*$
via the factorization map. Then ${\mathcal O}$ and ${\mathcal O}_-$ are
symplectic leaves in $G$ and $B_-$ respectively.

The Lie group $B_+$ acts on ${\mathcal O}$
as a subgroup of $G$. This action and the dressing action on
${\mathcal O}_-$ are quasi-Hamiltonian, in a sense that there is
an appropriate moment map \cite{Lu}.
Thus, the Lie group $B_+$ acts (locally) on
${\mathcal O}\times {\mathcal O}_-$
via the diagonal action and this action is quasi-Hamiltonian.
Therefore, we can reduce this product via Hamiltonian reduction and
thus we obtain the symplectic variety
\[
X({\mathcal O},{\mathcal O}_-)={\mathcal O}\times {\mathcal O}_-///B_+
\]
Here two dashes mean that we take the categorical quotient
and one extra dash means that we do Hamiltonian reduction.

It is easy to see that
${\rm dim}(X({\mathcal O},{\mathcal O}_-))=|\Delta_+|+{\rm rk}(w_0-id)-2n$.
Therefore the multiplicity of the restriction is
$\ell^{\frac{{\rm dim}(X({\mathcal O},{\mathcal O}_-))}{2}}$.

Comparing with  the above argument we arrive to a
conjecture:
\begin{conjecture} There is a birational correspondence between
the variety $({\mathcal O}\times {\mathcal O}_-)///B_+$
described above and a symplectic torus of the same dimension.
\end{conjecture}

{\bf 2}. The dimension in the tensor product. Consider three
conjugation $G$-orbits
${\mathcal O}_1,{\mathcal O}_2,{\mathcal O}_3\subset G$.
Consider $G$ as a Poisson variety with the Poisson structure
inherited from the dual Poisson Lie group $G^*$ via the factorization
map. This Poisson Lie structure is not a Poisson Lie structure
but rather a nonlinear deformation of the Poisson structure on
${\mathfrak g}^*$. Then conjugation orbits are
symplectic leaves and the natural action of $G$ on them is
quasi-Hamiltonian (i.e. there is an appropriate moment map).
This action induces a $G$-action on the product
${\mathcal O}_1\times {\mathcal O}_2\times \bar{{\mathcal O}}_3$
and this action is quasi-Hamiltonian. Here $\bar{{\mathcal O}}$
is the opposite symplectic variety to ${\mathcal O}$.

Consider the Hamiltonian reduction
\[
X({\mathcal O}_1,{\mathcal O}_2, {\mathcal O}_3)=({\mathcal
O}_1\times{\mathcal O}_2\times\bar{{\mathcal O}}_3)///G
\]
It is a symplectic variety with
${\rm dim}(X({\mathcal O}_1,{\mathcal O}_2,{\mathcal O}_3))=2(|\Delta_+|-n)$.

It is clear that the multiplicity of a generic irreducible
module in the tensor product of two generic irreducible
representations is
$\ell^{{\rm dim}(X({\mathcal O}_1,{\mathcal O}_2,{\mathcal O}_3))/2}$.

\begin{conjecture} There is a birational equivalence between the
symplectic variety
$({\mathcal O}_1\times{\mathcal O}_2\times\bar{{\mathcal O}}_3)///G$
and a complex torus of the same dimension.
\end{conjecture}

\end{document}